\documentclass[article]{IEEEtran}
\usepackage{cite}
\usepackage{color}
\usepackage[pdftex]{graphicx}
\graphicspath{{fig/}{jpeg/}}
\usepackage[cmex10]{amsmath}
\usepackage{amssymb}
\usepackage{algorithm}
\usepackage{algorithmic}
\usepackage{amsmath,amssymb,lipsum}
\usepackage{hyperref}
\usepackage{comment}
\usepackage{graphicx}
\usepackage[font=small]{caption}
\usepackage{subcaption}
\usepackage{setspace}

\input{mysymbol.sty}
\usepackage{needspace}

\newcommand{\myparagraph}[1]{\needspace{1\baselineskip}\medskip\noindent {\it #1.}}

\renewcommand{\QED}{\hfill\ensuremath{\blacksquare}}
\usepackage{theorem}
\newtheorem{thm}{Theorem}
\newtheorem{lemma}{Lemma}
\newtheorem{proposition}{Proposition}
\theoremstyle{definition}
\newtheorem{assumption}{Assumption}
\newtheorem{remark}{Remark}

\title{Network Newton--Part II: \\ Convergence Rate and Implementation}

\author{Aryan Mokhtari, Qing Ling and Alejandro Ribeiro 
\thanks{{Work supported by NSF CAREER CCF-0952867 and ONR N00014-12-1-0997. A. Mokhtari and A. Ribeiro are with the Department of Electrical and Systems Engineering, University of Pennsylvania, 200 South 33rd Street, Philadelphia, PA 19104. Q. Ling is with the Department of Automation, University of Science and Technology of China,
Hefei, Anhui 230026, China.  Email at: \{aryanm, aribeiro\}@seas.upenn.edu,  qingling@mail.ustc.edu.cn. Part of the results in this paper appeared in \cite{NNAsilomar} and \cite{NNICASSP}. {This paper expands the results and presents convergence proofs that are referenced in \cite{NNAsilomar} and \cite{NNICASSP}.} 
}}
}

\begin{document}

\maketitle
\thispagestyle{empty}

\begin{abstract}

The use of network Newton methods for the decentralized optimization of a sum cost distributed through agents of a network is considered. Network Newton methods reinterpret distributed gradient descent as a penalty method, observe that the corresponding Hessian is sparse, and approximate the Newton step by truncating a Taylor expansion of the inverse Hessian. Truncating the series at $K$ terms yields the NN-$K$ that requires aggregating information from $K$ hops away. Network Newton is introduced and shown to converge to the solution of the penalized objective function at a rate that is at least linear in a companion paper \cite{NN-part1}. The contributions of this work are: (i) To complement the convergence analysis by studying the methods' rate of convergence. (ii) To introduce adaptive formulations that converge to the optimal argument of the original objective. (iii) To perform numerical evaluations of NN-$K$ methods. The convergence analysis relates the behavior of NN-$K$ with the behavior of (regular) Newton's method and shows that the method goes through a quadratic convergence phase in a specific interval. The length of this quadratic phase grows with $K$ and can be made arbitrarily large. The numerical experiments corroborate reductions in the number of iterations and the communication cost that are necessary to achieve convergence relative to distributed gradient descent. 
\end{abstract}



\section{Introduction}\label{sec_Introduction}

In decentralized optimization problems a group of agents is tasked with minimizing a sum cost when each of them has access to a specific summand. They do so by working through subsequent rounds of local processing and variable exchanges with adjacent peers. This architecture arises naturally in decentralized control \cite{Bullo2009,Cao2013-TII,LopesEtal8} as well as in wireless \cite{Ribeiro10,Ribeiro12} and sensor networks \cite{Schizas2008-1,KhanEtal10,cRabbatNowak04}. In these problems agents have access to local information but want to achieve a common goal, administer a shared resource, or estimate the state of a global environment. Decentralized optimization is also relevant to large scale machine learning \cite{bekkerman2011scaling,Tsianos2012-allerton-consensus,Cevher2014}, where problems are not inherently distributed but are divvied up to process big datasets.

Irrespectively of the specific application, various methods have been developed for decentralized optimization. These include distributed gradient descent (DGD) \cite{Nedic2009,Jakovetic2014-1,YuanQing, Shi2014} as well as distributed implementations of the alternating direction method of multipliers \cite{Schizas2008-1,cQingRibeiroADMM14,BoydEtalADMM11,Shi2014-ADMM} and dual averaging \cite{Duchi2012,cTsianosEtal12}. At the core of all of these methods lies a gradient descent iteration that endows them with their convergence properties, but also results in large convergence times for problems with poor conditioning. In a companion paper, we introduced the network Newton family of decentralized optimization methods that incorporates second order information into DGD iterations to accelerate convergence \cite{NN-part1}. Methods in the network Newton family are derived by introducing a penalty formulation of distributed optimization objectives (Section \ref{sec:problem}) for which the resulting Hessians have the same sparsity pattern of the underlying network (Section \ref{subsec_nn}). The Hessian inverse that is necessary to compute Newton steps is then expressed as a Taylor series expansion that we truncate at $K$ terms to obtain the $K$th member of the network Newton family -- which we abbreviate as NN-$K$. These truncations can be computed in a distributed manner by aggregating information from, at most, $K$ hops away.

The network Newton methods have been proven to converge to the optimal solution of the penalized objective at a rate that is at least linear \cite{NN-part1}. The main goal of this paper is to complete the convergence analysis of NN-$K$ by studying its rate of convergence (Section \ref{sec:convergence_analysis}). We show that for all iterations except the first few, a weighted gradient norm associated with NN-$K$ iterates follows a decreasing path akin to the path that would be followed by regular Newton iterates (Lemma \ref{two_phase_convergence_lemma}). The only difference between these residual paths is that the NN-$K$ path contains a term that captures the error of the Hessian inverse approximation. Leveraging this similarity, it is possible to show that the rate of convergence is quadratic in a specific interval whose length depends on the order $K$ of the selected network Newton method (Theorem \ref{quadratic_convergence_theorem}). Existence of this quadratic convergence phase explains why NN-$K$ methods converge faster than DGD -- as we indeed observe in numerical analyses. It is also worth remarking that the error in the Hessian inverse approximation can be made arbitrarily small by increasing the method's order $K$ and, as a consequence, the quadratic phase can be made arbitrarily large.

Given that NN-$K$ solves the minimization of a penalized objective, it converges to a point that is close to the optimum. To achieve exact convergence we introduce an adaptive version of NN-$K$ -- which we term ANN-$K$ -- that uses a sequence of increasing penalty coefficients to achieve exact convergence to the optimal solution (Section \ref{sec:implement}). We wrap up the paper with numerical analyses. We first demonstrate the advantages of NN-$K$ relative to DGD for the minimization of a family of quadratic objective functions with varying condition number and network connectivity (Section \ref{sec:simulations}). As expected, NN-$K$ methods reduce convergence times by substantive factors when the objective functions are not well conditioned. Advantages in terms of communication cost are less marked because NN-$K$ aggregates information from $K$-hop neighborhoods but still substantial. Network Newton is also applied to solve a logistic regression problem. The results reinforce the conclusions reached for the quadratic objective problem (Section \ref{sec_logistic_regression}). Numerical analyses also show that network Newton methods with $K=1$ and $K=2$ tend to work best when measured in terms of overall communication cost. Numerical experiments for ANN-$K$ illustrate the tradeoffs that appear in the selection of the initial penalty coefficient and its rate of change (Section \ref{sec_sims_ann}). The paper closes with concluding remarks (Section \ref{sec_conclusions}).

\myparagraph{\bf Notation} Vectors are written as $\bbx\in\reals^n$ and matrices as $\bbA\in\reals^{n\times n}$. Given $n$ vectors $\bbx_i$, the vector $\bby=[\bbx_1;\ldots;\bbx_n]$ represents a stacking of the elements of each individual $\bbx_i$. The null space of matrix $\bbA$ is denoted by $\text{null}(\bbA)$ and the span of a vector by $\text{span}(\bbx)$. The $i$-th eigenvalue of matrix $\bbA$ is denoted by $\mu_{i}(\bbA)$. For matrices $\bbA$ and $\bbB$ their Kronecker product is denoted as $\bbA\otimes\bbB$. The gradient of a function $f(\bbx)$ is denoted as $\nabla f(\bbx)$ and the Hessian matrix is denoted by $\nabla^2 f(\bbx)$.


%
\section{Algorithm definition} \label{sec:problem}

We consider a connected and symmetric network with $n$ agents generically indexed by $i=1,\ldots,n$. The network is specified by the $n$ neighborhood sets $\mathcal{N}_i$, each of which is defined as the group of nodes that are connected to $i$. Nodes have access to strongly convex local objective functions $f_i(\bbx)$, but cooperate to minimize the global cost $f(\bbx):=\sum_{i=1}^{n}\ f_{i}(\bbx)$, 
\begin{equation}\label{original_optimization_problem1}
  \bbx^* \ :=\ \argmin_{\bbx}  f(\bbx) 
      \  =\ \argmin_{\bbx}\sum_{i=1}^{n} f_i(\bbx).
\end{equation}
To rewrite the global problem in a form that is suitable for distributed implementation we define local variables $\bbx_i\in \reals^p$ and rewrite the cost to be minimized as $\sum_{i=1}^{n}\ f_{i}(\bbx_i)$. For a problem formulation equivalent to \eqref{original_optimization_problem1}, we have to further add the restriction that local variables $\bbx_{i}$ be the same as neighboring variables $\bbx_j$ with $j\in\ccalN_i$,
\begin{align}\label{original_optimization_problem2}
   \{\bbx_i^*\}_{i=1}^n\ := \
   &\argmin_{\bbx} \ \sum_{i=1}^{n}\ f_{i}(\bbx_{i}), \nonumber\\ 
   &\text{\ s.t.}  \ \bbx_{i}=\bbx_{j}, 
                   \quad \text{for all\ } i, j\in\ccalN_i .
\end{align} 
The optimization problems in \eqref{original_optimization_problem2} and \eqref{original_optimization_problem1} are equivalent in the sense that $\bbx_i^*=\bbx^*$ for all $i$. This has to be true because the feasible set of \eqref{original_optimization_problem2} is restricted to configurations in which all variables $\bbx_i$ are equal given that the network is connected.

The constraints $\bbx_{i}=\bbx_{j}$ imposed for all $i$ and all $j\in\ccalN_i$ are a way of making all local variables equal but there are other alternatives. The one that is germane to this paper consists of introducing weights $w_{ij}$ that we group in the matrix $\bbW\in\reals^{n\times n}$. The weights $w_{ij}$ are chosen so that the matrix $\bbW$ is symmetric, row stochastic, and such that the null space of $\bbI-\bbW$ is the span of the all one vector $\bbone$
\begin{equation}\label{eqn_conditions_on_weights}
   \bbW^T=\bbW, \quad 
   \bbW\bbone=\bbone, \quad
   \text{null}(\bbI-\bbW)=\text{span}(\bbone).
\end{equation}
We further define the extended weight matrix 
\begin{equation}\label{eqn_extended_weight_matrix}
   \bbZ:= \bbW \otimes \bbI\in\reals^{np\times np}
\end{equation}
as the Kronecker product of the weight matrix $\bbW\in\reals^{n\times n}$ and the identity matrix $\bbI\in\reals^{p\times p}$ as well as the vector $\bby := \left[\bbx_{1}; \dots ; \bbx_{n}\right]$ as the concatenation of the local vectors $\bbx_{i}$. It follows that the equality constraint $\bbZ= \bbW \otimes \bbI$ can be satisfied if and only if all the local variables are equal, i.e., if and only if $\bbx_1=\dots=\bbx_n$. Indeed, since the null space of $\bbI-\bbW$ is $\text{null}(\bbI-\bbW)=\text{span}(\bbone)$ as per the last condition in \eqref{eqn_conditions_on_weights}, the null space of  $\bbI-\bbZ$ must be $\text{null}(\bbI-\bbZ)=\text{span}(\bbone\otimes\bbI)$. Thus, vectors $\bby := \left[\bbx_{1}; \dots ; \bbx_{n}\right]$ in the null space of $\bbI-\bbZ$, which, by definition, are the only ones that satisfy the equality $(\bbI-\bbZ)\bby=\bb0$ are multiples of $\bbone\otimes\bbI$ and therefore satisfy $\bbx_1=\dots=\bbx_n$. 

Further observe that the matrix $\bbZ$, being stochastic and symmetric, is positive semidefinite. Consequently, the square root matrix $({\bbI-\bbZ})^{1/2}$ exists and has the same null space of $\bbI-\bbZ$. It then follows that $({\bbI-\bbZ})^{1/2}\bby=\bbzero$ if and only if the components of $\bby$ satisfy $\bbx_1=\dots=\bbx_n$. In turn, this implies that the optimization problem in \eqref{original_optimization_problem2} is equivalent to
\begin{align}\label{original_optimization_new_notation}
   \tby^*\ :=\ &\argmin_{\bbx}\ \sum_{i=1}^{n}\ f_i(\bbx_i), \nonumber \\
               &\text{\ s.t.} \quad ({\bbI-\bbZ})^{1/2}  \bby =\bb0.
\end{align} 
Here, we solve \eqref{original_optimization_new_notation} using a penalized version of the objective function. To do so we consider a given penalty coefficient $1/\alpha$ and the squared norm penalty function $(1/2)\|({\bbI-\bbZ})^{1/2}  \bby \|^2 = (1/2)\bby^{T}(\bbI -\bbZ) \bby$ associated with the constraint $({\bbI-\bbZ})^{1/2}  \bby =\bb0$. With penalty function and coefficient so defined, we can now introduce the penalized objective $F(\bby):=(1/2)\bby^{T}(\bbI -\bbZ) \bby + \alpha\sum_{i=1}^{n} f_i(\bbx_i)$ and the penalized optimization problem 
\begin{align}\label{centralized_opt_problem}
    \bby^* :=& \argmin\ F(\bby) \nonumber\\
           :=& \argmin \frac{1}{2}\ \bby^{T}(\bbI -\bbZ)\ \bby 
                       + \alpha\sum_{i=1}^{n} f_i(\bbx_i).
\end{align}
As the penalty coefficient $1/\alpha$ grows, or, equivalently, as $\alpha$ vanishes, the optimal argument $\bby^*$ of the penalized problem \eqref{centralized_opt_problem} converges towards the optimal argument $\tby^*$ of \eqref{original_optimization_problem2} and \eqref{original_optimization_new_notation}. In that sense, \eqref{centralized_opt_problem} is a reasonable proxy for \eqref{original_optimization_problem2}, \eqref{original_optimization_new_notation}, and the equivalent original formulation in \eqref{original_optimization_problem1}.

The property that makes the penalized problem in \eqref{centralized_opt_problem} amenable to distributed implementation is that its gradients can be computed by exchanging information between neighboring nodes. This property is the basis for the development of the DGD method of \cite{Nedic2009} and the NN method of \cite{NN-part1}. In the following section we study the idea of using Newton's method for solving \eqref{centralized_opt_problem}.

%
\subsection{Newton's method and Hessian splitting}\label{subsec_nn}

We proceed to minimize the penalized objective function $F(\bby)$ in \eqref{centralized_opt_problem} using Newton's method. The Newton update with stepsize $\epsilon $ for function $F(\bby)$ can be written as
\begin{equation}\label{newton_step}
\bby_{t+1}=\bby_t-\eps\ \! \nabla^2F(\bby_t)^{-1}   \nabla F(\bby_t),
\end{equation}
where $\nabla^2F(\bby_t)$ and $\nabla F(\bby_t)$ are the Hessian and gradient of function $F$ evaluated at point $\bby_t$, respectively. 

To compute the gradient $\nabla F(\bby_t)$ we introduce the vector $\bbh(\bby):=[\nabla f_{1}(\bbx_1);\dots;\nabla f_{n}(\bbx_n)]$ that concatenates the local gradients $\nabla f_{i}(\bbx_i)$. Given the definition of the penalized function $F(\bby)$ in \eqref{centralized_opt_problem} it follows that the gradient of $F(\bby)$ at $\bby=\bby_t$ is
\begin{equation}\label{eqn_gradient_definition}
    \bbg_t \ :=\ \nabla F(\bby_t) 
           \  =\ (\bbI -\bbZ) \bby_t +\alpha \bbh(\bby_t).
\end{equation}
The computation of the gradient $\bbg_t$ can be distributed through the network because $\bbZ$ has the sparsity pattern of the graph. Specifically, define the local gradient component $\bbg_{i,t}$ as the $i$th element of the gradient $\bbg_{t}=[\bbg_{i,t};\ldots;\bbg_{i,t}]$ and recall that $\bbx_{i,t}$ and $\bbx_{i,t+1}$ are the $i$th components of the vector $\bby_t$ and $\bby_{t+1}$. The local gradient component at node $i$ is given by
\begin{equation}\label{local_gradient}
\bbg_{i,t}=(1-w_{ii})\bbx_{i,t} - \sum_{j\in \mathcal{N}_i} w_{ij} \bbx_{j,t}+\alpha \nabla f_{i}(\bbx_{i,t}).
\end{equation}
Using \eqref{local_gradient}, node $i$ can compute its local gradient using its local iterate $\bbx_{i,t}$, the gradient of the local function $\nabla f_{i}(\bbx_{i,t})$ and the $\bbx_{j,t}$ iterates of its neighbors $j\in \ccalN_i$.

To implement Newton's method as defined in \eqref{newton_step} we also need to compute the Hessian $\bbH_t:=\nabla^2 F(\bby_t)$ of the penalized objective. Start by differentiating twice the objective function $F$ in \eqref{centralized_opt_problem} in order to write the Hessian $\bbH_t$ as
\begin{equation}\label{Hessian}
   \bbH_t := \nabla^2 F(\bby_t) = \bbI-\bbZ +\alpha \bbG_t,
\end{equation}
where the matrix $\bbG_t \in \reals^{np\times np}$ is a block diagonal matrix formed by blocks $\bbG_{ii,t}\in \reals^{p\times p}$ containing the Hessian of the $i$th local function,
\begin{equation}\label{G_form}
   \bbG_{ii,t} = \nabla^2 f_i(\bbx_{i,t})  .
\end{equation}
It follows from \eqref{Hessian} and \eqref{G_form} that the Hessian $\bbH_t$ is block sparse with blocks $\bbH_{ij,t}\in \reals^{p\times p}$ having the sparsity pattern of $\bbZ$, which is the sparsity pattern of the graph. {The diagonal blocks are of the form $\bbH_{ii,t}=(1-w_{ii})\bbI +  \alpha \nabla^2 f_i(\bbx_{i,t}) $ and the off diagonal blocks are not null only when $j\in \mathcal{N}_i$ in which case $\bbH_{ij,t}=w_{ij}\bbI$.}

Recall that for the Newton update in \eqref{newton_step}, the Hessian inverse $\nabla^2F(\bby_t)^{-1}=\bbH_t^{-1}$ evaluated at $\bby=\bby_t$ is required not the Hessian $\bbH_t$. While the Hessian $\bbH_t$ is sparse, the inverse $\bbH_t^{-1}$ is not necessarily sparse. Therefore, the Hessian inverse is not necessarily computable in a decentralized setting. To overcome this problem we split the diagonal and off diagonal blocks of $\bbH_t$ and rely on the Taylor's expansion of the inverse $\bbH_t^{-1}$. To be precise, write $\bbH_t=\bbD_t - \bbB$ where the matrix $\bbD_t$ is defined as 
\begin{equation}\label{diagonal_matrix}
   \bbD_t := \alpha \bbG_t + 2\ ( \bbI  -  \diag(\bbZ))
          := \alpha \bbG_t + 2\ ( \bbI  -  \bbZ_{d}).
\end{equation}
In the second equality we defined $\bbZ_{d}:=\diag(\bbZ)$ for future reference. Observe that the matrix $\bbI-\bbZ_{d}$ is positive definite because in a connected network the local weights are $w_{ii}<1$. The block diagonal matrix $\bbG_t$ is also positive definite because the local functions are assumed strongly convex. it follows from these two observations that the matrix $\bbD_t$ is block diagonal and positive definite. Further note that the $i$th diagonal block $\bbD_{ii,t}\in\reals^p$ of $\bbD_t$ can be computed and stored by node $i$ as $\bbD_{ii,t}= \alpha \nabla^2 f_{i}(\bbx_{i,t}) + 2(1-w_{ii})\bbI $ using local information only. To have $\bbH_t=\bbD_t - \bbB$ we must define $\bbB:=\bbD_t-\bbH_t$.  Considering the definitions of $\bbH_t$ and $\bbD_t$ in \eqref{Hessian} and \eqref{diagonal_matrix}, respectively, it follows that 
\begin{equation}\label{non_diagona_matrix}
   \bbB =  \bbI - 2\bbZ_{d} +\bbZ.
\end{equation}
Observe that $\bbB$ is independent of time and only depends on the weight matrix $\bbZ$. As in the case of the Hessian $\bbH_t$, the matrix $\bbB$ is block sparse with blocks $\bbB_{ij}\in \reals^{p\times p}$ having the sparsity pattern of $\bbZ$, which is the sparsity pattern of the graph. Since $\bbB$ is block sparse, node $i$ can compute the diagonal block $\bbB_{ii}=(1-w_{ii})\bbI$ and the off diagonal blocks $\bbB_{ij}=w_{ij}\bbI$ using local information about its own weights.

Proceed now to factor $\bbD_t^{1/2}$ from both sides of the splitting relationship to write $\bbH_t = \bbD_t ^{{1}/{2}} (  \bbI - \bbD_t ^{-{1}/{2}}\bbB\bbD_t ^{-{1}/{2}} ) \bbD_t^{{1}/{2}}$. This decomposition implies that the Hessian inverse $\bbH_t^{-1}$ can be computed from the Taylor series expansion $(\bbI-\bbX)^{-1}= \sum_{j=0}^{\infty} \bbX^{j}$ with $\bbX=\bbD_t^{-{1}/{2}}  \bbB   \bbD_t^{-{1}/{2}}$. Therefore, we can write 
\begin{equation}\label{exact_Hessian_inverse}
   \bbH_t^{-1} = \bbD_t^{-1/2} 
                \sum_{k=0}^{\infty} \left(\bbD_t^{-1/2}  
                  \bbB   \bbD_t^{-1/2}\right)^{k}\ \bbD_t^{-1/2}.
\end{equation}
The sum in \eqref{exact_Hessian_inverse} converges if the absolute value of all the eigenvalues of the matrix $\bbD^{-{1}/{2}}  \bbB   \bbD^{-{1}/{2}} $ are strictly less  than 1 -- we prove that this is true in Proposition \ref{symmetric_term_bounds11}. Truncations of this convergent series are utilized to define the family of Network Newton methods in the following section.

%
\subsection{Network Newton}\label{sec_the_method}

Network Newton is defined as a family of algorithms that rely on truncations of the series in \eqref{exact_Hessian_inverse}. The {$K$th} member of this family, NN-$K$, considers the first $K+1$ terms of the series to define the approximate Hessian inverse
\begin{equation}\label{Hessian_inverse_approximation}
   \hbH_t^{(K)^{-1}} :=     \bbD_t^{-1/2}  \  \sum_{k=0}^{K} \left(  \bbD_t^{-1/2}  \bbB   \bbD_t^{-1/2}            \right)^{k}     \ \bbD_t^{-1/2}.
\end{equation}
NN-$K$ uses the approximate Hessian $\hbH_t^{(K)^{-1}}$ as a curvature  correction matrix that is used in lieu of the exact Hessian inverse $\bbH^{-1}$ to estimate the Newton step. I.e., instead of descending along the Newton step $\bbd_t:=-\bbH_t^{-1}\bbg_t$ we descend along the NN-$K$ step $\bbd_t^{(K)}:=-\hbH_t^{(K)^{-1}}\bbg_t$, which we intend as an approximation of $\bbd_t$. Using the explicit expression for $\hbH_t^{(K)^{-1}}$ in \eqref{Hessian_inverse_approximation} we write the NN-$K$ step as
\begin{equation}\label{Hessian_approximation_iteration}
\bbd_t^{(K)} = -\  \bbD_t^{-1/2}  \  \sum_{k=0}^{K} \left(  \bbD_t^{-1/2}  \bbB   \bbD_t^{-1/2}            \right)^{k}     \ \bbD_t^{-1/2}\ \bbg_t,
\end{equation}
where, we recall, the vector $\bbg_t$ is the gradient of the objective function $F(\bby)$ defined in \eqref{eqn_gradient_definition}. The NN-$K$ update formula can then be written as
\begin{equation}\label{update_formula_NN}
   \bby_{t+1}=\bby_t+\eps\  \bbd_t^{(K)}.
\end{equation}
The algorithm defined by recursive application of \eqref{update_formula_NN} can be implemented in a distributed manner. Specifically, define the components $\bbd^{(K)}_{i,t}\in\reals^p$ of the NN-$K$ step $\bbd^{(K)}_{t}=[\bbd^{(K)}_{1,t};\ldots;\bbd^{(K)}_{n,t}]$ and rewrite \eqref{update_formula_NN} componentwise as 
\begin{equation}\label{update_formula_NN_local}
   \bbx_{i,t+1}=\bbx_{i,t} +\eps\  \bbd^{(K)}_{i,t}.
\end{equation}
To determine the step components $\bbd^{(K)}_{i,t}$ in \eqref{update_formula_NN_local} observe that considering the definition of the NN-$K$ descent direction in \eqref{Hessian_approximation_iteration}, Network Newton descent directions can be computed by the recursive expression 
\begin{equation}\label{eqn_nn_step_recursion}
   \bbd_t^{(k+1)} = \bbD_t^{-1}\bbB \bbd_t^{(k)} -\bbD_t^{-1}\bbg_t
                  = \bbD_t^{-1}\left(\bbB \bbd_t^{(k)} - \bbg_t \right).
\end{equation}
%
%
If we expand the product $\bbB \bbd_t^{(k)}$ as a sum and utilize the fact that the blocks of the matrix $\bbB$ have the sparsity pattern of the graph, we can separate \eqref{eqn_nn_step_recursion} into the following componentwise recursions
\begin{equation}\label{local_descent}
  \bbd_{i,t}^{(k+1)} 
      = \bbD_{ii,t}^{-1}\bigg[
              \sum_{j\in \mathcal{N}_i,j=i} \bbB_{ij} \bbd_{j,t}^{(k)} 
              - \bbg_{i,t}\bigg].
\end{equation}
That the matrix $\bbB$ is block-sparse permits writing the sum in \eqref{local_descent} as a sum over neighbors, instead of a sum across all nodes. 

In \eqref{local_descent}, the matrix blocks $\bbD_{ii,t}=\alpha \nabla^2 f_{i}(\bbx_{i,t}) + 2(1-w_{ii})\bbI$, $\bbB_{ii}=(1-w_{ii})\bbI$, and $\bbB_{ij}=w_{ij}\bbI$ are evaluated and stored at node $i$. The gradient component $\bbg_{i,t}$ is also stored and computed at $i$ upon being communicated the values of neighboring iterates $\bbx_{j,t}$ [cf. \eqref{local_gradient}]. Thus, if the NN-$k$ step components $\bbd_{j,t}^{(k)}$ are available at neighboring nodes $j$, node $i$ can determine the  NN-$(k+1)$ step component $\bbd_{i,t}^{(k+1)}$ upon being communicated that information. We use this property to embed an iterative computation of the NN-$K$ step inside the NN-$K$ recursion in \eqref{update_formula_NN_local}. For each iteration index $t$, we compute the local component of the NN-$0$ step $\bbd_{i,t}^{(0)}=-\bbD_{ii,t}^{-1}\bbg_{i,t}$. Upon exchanging this information with neighbors we use \eqref{local_descent} to determine the NN-$1$ step components $\bbd_{i,t}^{(1)}$. These can be exchanged and plugged in \eqref{local_descent} to compute $\bbd_{i,t}^{(2)}$. Repeating this procedure $K$ times, nodes end up having determined their NN-$K$ step component $\bbd_{i,t}^{(K)}$. They use this step to update $\bbx_{i,t}$ according to \eqref{update_formula_NN_local} and move to the next iteration. We analyze the convergence rate of the resulting algorithm in Section \ref{sec:convergence_analysis} and develop a numerical analysis in Section \ref{sec:simulations}.

%
\begin{remark} \label{rmk_comm_cost}\normalfont
By trying to approximate the Newton step, NN-$K$ ends up reducing the number of iterations required for convergence. Furthermore, the larger $K$ is, the closer that the NN-$K$ step gets to the Newton step, and the faster NN-$K$ converges. We will justify these assertions both, analytically in Section \ref{sec:convergence_analysis}, and numerically in Section \ref{sec:simulations}. It is important to observe, however, that reducing the number of iterations reduces the computational cost but not necessarily the communication cost. In DGD, each node $i$ shares its vector $\bbx_{i,t}\in \reals^{p}$ with each of its neighbors $j \in \mathcal{N}_i$. In NN-$K$, node $i$  exchanges not only the vector $\bbx_{i,t}\in \reals^{p}$ with its neighboring nodes, but it also communicates iteratively the local components of the descent directions $\{\bbd_{i,t}^{(k)}\}_{k=0}^{K-1}\in \reals^{p}$ so as to compute the descent direction $\bbd_{i,t}^{(K)}$. Therefore, at each iteration, node $i$ sends $|\mathcal{N}_{i}|$ vectors of size $p$ to the neighboring nodes in DGD, while in NN-$K$ it sends $(K+1)|\mathcal{N}_{i}|$ vectors of the same size. Unless the original problem is well conditioned, NN-$K$ also reduces total communication cost until convergence, even though the cost of each individual iteration is larger. However, the use of large $K$ is unwarranted because the added benefit of better approximating the Newton step does not compensate the increase in communication cost. 
\end{remark}

%
\section{Convergence Rate}\label{sec:convergence_analysis}

Linear convergence of the sub optimality sequence $F(\bby_t)-F(\bby^*)$ associated with NN-$K$ iterates $\bby_t$ has been proven in \cite{NN-part1}. We improve this result by showing that regardless of the choice of $K$, the rate of convergence is quadratic in a specific interval. To prove this result we utilize some of the results in \cite{NN-part1} which we repeat here for completeness. We start by restating assumptions that are necessary for the convergence analysis.

%
\begin{assumption}\label{ass_weight_bounds} There exists constants $0\leq\delta\leq\Delta<1$ that lower and upper bound the diagonal weights for all $i$, 
\begin{equation}\label{bounds_for_local_weights}
   0\leq \delta  \leq w_{ii} \leq \Delta <1,  \qquad  i=1,\ldots,n .
\end{equation}\end{assumption}

%
\begin{assumption}\label{convexity_assumption} 
The local objective functions $f_i(\bbx)$ are twice differentiable and the eigenvalues of the local objective function Hessians are bounded with positive constants $0<m\leq M<\infty$, i.e. 
\begin{equation}\label{local_hessian_eigenvlaue_bounds}
m\bbI\preceq \nabla^2 f_i(\bbx)\preceq M\bbI.
\end{equation}
\end{assumption}

%
\begin{assumption}\label{Lipschitz_assumption} The local objective function Hessians $\nabla^2 f_i(\bbx)$ are Lipschitz continuous with respect to the Euclidian norm with parameter $L$. I.e., for all $\bbx, \hbx \in \reals^p$, it holds
\begin{equation}
   \left\| \nabla^2 f_i(\bbx)-\nabla^2 f_i(\hbx) \right\| \ \leq\  L\ \| \bbx- \hbx \|.
\end{equation}
\end{assumption}

%
The upper bound $\Delta<1$ on the local weights $w_{ii}$ in Assumption \ref{ass_weight_bounds} exits for connected networks. The non-negative lower bound  $\delta$ on the local weights $w_{ii}$ is more a definition than a constraint since we may have $\delta=0$. Strong convexity of the local objective functions $f_i$ enforces the existence of a lower bound $m$ for the eigenvalues of the local Hessian $\nabla^2f_i$ as in \eqref{local_hessian_eigenvlaue_bounds}. The upper bound $M$ for the eigenvalues of local objective function Hessians $\nabla^2 f_i(\bbx)$ in Assumption \ref{convexity_assumption} is equivalent to the assumption that local gradients $\nabla f_i(\bbx)$ are Lipschitz continuous with parameter $M$. Assumption \ref{Lipschitz_assumption} states that the local objective function Hessians are Lipschitz continuous with parameter $L$. A particular consequence of this assumption is that the penalized objective function Hessian $\bbH(\bby):= \nabla^2 F(\bby)$ is also Lipschitz continuous with parameter $\alpha L$ -- see Lemma 1 of \cite{NN-part1}. I.e. for all $\bby,\hby\in\reals^{np}$ it holds,
\begin{equation}\label{H_Lipschitz_claim}
    \left\|\bbH(\bby)-\bbH(\hby)\right\| \leq \alpha L \| \bby-\hby\| ,
\end{equation}
Recall that the block diagonal matrix $\bbD_{t}$, being the sum of positive definite $\alpha \bbG_{t}$ and $2( \bbI  -  \bbZ_{d})$, is positive definite and, therefore, invertible. Further recall that the matrix $\bbB$, being symmetric and doubly stochastic, has eigenvalues that lie between 0 and 1 and is therefore positive semidefinite. These facts can be used to prove that the eigenvalues of the matrix $\bbD_t^{-{1}/{2}}\bbB\bbD_t^{-{1}/{2}}$ must be nonnegative and strictly smaller than $1$ as we state next  \cite[Proposition 2]{NN-part1}.

%
\begin{proposition}\label{symmetric_term_bounds11}
Consider the definitions of matrices $\bbD_t$ in \eqref{diagonal_matrix} and $\bbB$ in \eqref{non_diagona_matrix}. If Assumptions \ref{ass_weight_bounds} and \ref{convexity_assumption} hold true, the matrix $\bbD_t^{-{1}/{2}}  \bbB   \bbD_t^{-{1}/{2}}$ is positive semidefinite and the eigenvalues are bounded above by a constant $\rho<1$
\begin{equation}\label{important_claim}
    \bb0\  \preceq \  \bbD_t^{-{1}/{2}}  \bbB   \bbD_t^{-{1}/{2}} \ \preceq\ \rho\bbI ,
\end{equation}
where $\rho:= 2(1-\delta)/(2(1-\delta)+{\alpha m} )$.
\end{proposition}

%
The result in Proposition \ref{symmetric_term_bounds11} makes the expansion in \eqref{exact_Hessian_inverse} valid and is used in subsequent proofs. These proofs also rely on guarantees that the eigenvalues of the approximate Hessian inverse $\hbH_t^{(K)^{-1}}$ are positive and finite for all choices of $K$ and for all steps $t$. We state these guarantees next \cite[Lemma 2]{NN-part1}.

%
\begin{lemma}\label{Hessian_inverse_eigenvalue_bounds_lemma}
Consider the NN-$K$ method as defined by \eqref{Hessian_approximation_iteration}-\eqref{update_formula_NN} with the gradient $\bbg_t$ as defined in \eqref{eqn_gradient_definition} and the matrices $\bbB$ and $\bbD_t$ defined as in \eqref{G_form}-\eqref{non_diagona_matrix}. If Assumptions \ref{ass_weight_bounds} and \ref{convexity_assumption} hold true, the eigenvalues of approximate Hessian inverse $\hbH_t^{(K)^{-1}}$ are bounded as
\begin{equation}\label{bounded_Hessian_inverse}
\lambda \bbI\ \preceq\  \hbH_t^{(K)^{-1}}  \preceq\   \Lambda \bbI,
\end{equation} 
where constants $\lambda$ and $\Lambda$ are defined as
\begin{equation}\label{definition_of_lambdas}
\!\!\! \lambda\!:=\! \frac{1}{2(1-\delta)+\alpha M } \ \text{ and} 
\ \ \Lambda\! :=\! {\frac{1-\rho^{K+1}}{(1-\rho)(2(1-\Delta)+\alpha m )}}.
\end{equation}
\end{lemma}

%
The lower bound $\lambda>0$ for the eigenvalues of the approximate Hessian inverse $ \hbH_t^{(K)^{-1}} $ guarantees decrement in each network Newton iteration. The upper bound $\Lambda<\infty$ ensures that the norm of the network Newton step $\|\hbH_t^{(K)^{-1}} \bbg_t\|$ is bounded by a factor proportional to the gradient norm $\|\bbg_t\|$. Both of these results are necessary to show that the network Newton direction $\hbH_t^{(K)^{-1}}\bbg_t$ is a descent direction. This is claimed to be true in the following theorem \cite[Theorem 1]{NN-part1}.

%
\begin{thm}\label{linear_convergence}
Consider the objective function $F(\bby)$ as introduced in \eqref{centralized_opt_problem} and the NN-$K$ method as defined by \eqref{Hessian_approximation_iteration}-\eqref{update_formula_NN} with the gradient $\bbg_t$ as defined in \eqref{eqn_gradient_definition} and the matrices $\bbB$ and $\bbD_t$ defined as in \eqref{G_form}-\eqref{non_diagona_matrix}. If the stepsize $\epsilon$ is chosen as 
\begin{equation}\label{step_size_condition}
\epsilon = \min  \left\{ 1\ , \left[{\frac{3m\lambda^{\frac{5}{2}}}{ L\Lambda^{3}{(F(\bby_{0})-F(\bby^{*})) }^{\frac{1}{2}} } }  \right]^{\frac{1}{2}} \right\} 
\end{equation}
and Assumptions \ref{ass_weight_bounds}, \ref{convexity_assumption}, and \ref{Lipschitz_assumption} hold true, the sequence $F(\bby_{t})$ converges to the optimal argument $F(\bby^{*})$ at least linearly with constant $1-\zeta$. I.e.,
\begin{equation}\label{linear_convegrence_claim}
F(\bby_{t}) -F(\bby^*) \leq (1-\zeta)^t  {\left(F(\bby_{0}) -F(\bby^*) \right)},
\end{equation}
where the constant $0<\zeta<1$ is explicitly given by
\begin{equation}\label{beta_0_defintion}
\zeta :=    {(2-\eps) \epsilon\alpha m\lambda} - \frac{\alpha\epsilon^3  L\Lambda^{3}(F(\bby_{0})-F(\bby^{*}) )^{\frac{1}{2}}}{6\lambda^{\frac{3}{2}}} .
\end{equation}\end{thm}

%
Theorem \ref{linear_convergence} establishes linear convergence of the sequence of penalized objective functions $F(\bby_t)$ generated by NN-$K$ to the optimal objective function $F(\bby^*)$ -- which implies convergence of $\bby_t$ to the optimal argument $\bby^*$. This result is identical to the convergence behavior of DGD as shown in, e.g., \cite{YuanQing}. We expect to observe faster convergence for NN-$K$ relative to DGD, since NN-$K$ uses an approximation of the curvature of the penalized objective function $F$. In the following section we show that this expectation is fulfilled and that NN-$K$ has a quadratic convergence phase regardless of the choice of $K$. 

%
\subsection{Quadratic convergence phase}

To characterize convergence rate of NN-$K$, we first study the difference between this algorithm and (exact) Newton's method. In particular, the following lemma shows that the convergence of the norm of the weighted gradient $\|\bbD_{t-1}^{-{1}/{2}}\bbg_{t}\|$ in NN-$K$ is akin to the convergence of Newton's method with constant stepsize. The difference is the appearance of a term associated with the error of the Hessian inverse approximation as we formally state next.

%
\begin{lemma}\label{two_phase_convergence_lemma}
Consider the NN-$K$ method as defined by \eqref{Hessian_approximation_iteration}-\eqref{update_formula_NN} with the gradient $\bbg_t$ as defined in \eqref{eqn_gradient_definition} and the matrices $\bbB$ and $\bbD_t$ defined as in \eqref{G_form}-\eqref{non_diagona_matrix}. If Assumptions \ref{ass_weight_bounds}, \ref{convexity_assumption}, and \ref{Lipschitz_assumption} hold true, the sequence of weighted gradients $\bbD_{t}^{-{1}/{2}}\bbg_{t+1}$ satisfies 
\begin{align}\label{convg_rate_lemma_claim}
\left\|\bbD_{t}^{-{1}/{2}}\bbg_{t+1}\right\|&   \leq \\
 &	   \!\!\!\! \!\!\!\!  \!\!\!\! \left(1-\eps+ \eps \rho^{K+1}\right) \left[1 + \Gamma_1 (1-\zeta)^{\frac{(t-1)}{4}} \right]\left\|\bbD_{t-1}^{-{1}/{2}}\bbg_{t}\right\| \nonumber \\	
	   &  + \eps^2 \Gamma_2 \left\| \bbD_{t-1}^{-{1}/{2}}\bbg_{t}\right\|^2,\nonumber
\end{align}
where the constants $\Gamma_1$ and $\Gamma_2$ are defined as 
\begin{align}\label{Gammas_definition}
   & \Gamma_1 := {\frac{(\alpha \eps L\Lambda)^{{1}/{2}}  (F(\bby_{0}) -F(\bby^*) )^{{1}/{4}}}{\lambda^{{3}/{4}}(2(1-\Delta)+\alpha m) }}, \nonumber \\
	&\Gamma_2 := \frac{ \alpha L\Lambda^2}{2\lambda {(2(1-\Delta)+\alpha m)}^{{1}/{2}}}.
\end{align}
\end{lemma}

%
\begin{myproof} See Appendix \ref{App_two_phase_convergence_lemma}. \end{myproof}

%
As per Lemma \ref{two_phase_convergence_lemma} the weighted gradient norm $\|\bbD_{t}^{-{1}/{2}}\bbg_{t+1}\|$ is upper bounded by terms that are linear and quadratic on the weighted norm $\|\bbD_{t-1}^{-{1}/{2}}\bbg_{t}\|$ associated with the previous iterate. This is akin to the gradient norm decrease of Newton's method with constant stepsize. To make this connection clearer, further note that for all except the first few iterations the term $\Gamma_1 (1-\zeta)^{(t-1)/4}\approx 0$ is close to 0 and the relation in \eqref{convg_rate_lemma_claim} can be simplified to
\begin{align}\label{convg_rate_lemma_claim_large_t}
\left\|\bbD_{t}^{-{1}/{2}}\bbg_{t+1}\right\|  \ \lesssim\ & \left(1-\eps+ \eps \rho^{K+1}\right)\left\|\bbD_{t-1}^{-{1}/{2}}\bbg_{t}\right\| \nonumber \\	
	   & \quad  + \eps^2 \Gamma_2 \left\| \bbD_{t-1}^{-{1}/{2}}\bbg_{t}\right\|^2.
\end{align}
In \eqref{convg_rate_lemma_claim_large_t}, the coefficient in the linear term is reduced to $(1-\eps+ \eps \rho^{K+1})$ and the coefficient in the quadratic term stays at $\eps^2 \Gamma_2$. If, for discussion purposes, we set $\eps=1$ as in Newton's quadratic phase, we see that the upper bound in \eqref{convg_rate_lemma_claim_large_t} is further reduced to 
\begin{align}\label{convg_rate_lemma_claim_large_t_eps_1}
\left\|\bbD_{t}^{-{1}/{2}}\bbg_{t+1}\right\|  \ \lesssim\ & \rho^{K+1}\| \bbD_{t-1}^{-{1}/{2}}\bbg_{t}\|+ \Gamma_2\| \bbD_{t-1}^{-{1}/{2}}\bbg_{t}\|^2
\end{align}
We do not obtain quadratic convergence as in Newton's method because of the term $\rho^{K+1}\| \bbD_{t-1}^{-{1}/{2}}\bbg_{t}\|$. However, since the constant $\rho$ (cf. Proposition \ref{symmetric_term_bounds11}) is smaller than 1 the term $\rho^{K+1}$ can be made arbitrarily small by increasing the approximation order $K$. Equivalently, this means that by selecting $K$ to be large enough, we can make the quadratic term in \eqref{convg_rate_lemma_claim_large_t_eps_1} dominant and observe a quadratic convergence phase. The boundaries of this quadratic convergence phase are formally determined in the following Theorem.

%
\begin{thm}\label{quadratic_convergence_theorem}
Consider the NN-$K$ method as defined by \eqref{Hessian_approximation_iteration}-\eqref{update_formula_NN} with the gradient $\bbg_t$ as defined in \eqref{eqn_gradient_definition} and the matrices $\bbB$ and $\bbD_t$ defined as in \eqref{G_form}-\eqref{non_diagona_matrix}. Define the sequence $\eta_{t} :=[(1-\eps+\eps\rho^{K+1})(1 +  \Gamma_1 (1-\zeta)^{{(t-1)}/{4}} ) ]$ and the time $t_0$ as the first time at which sequence $\eta_t$ is smaller than 1, i.e. $t_0:= \argmin_{t}\{ t \mid \eta_t<1\}$. If Assumptions  \ref{ass_weight_bounds}, \ref{convexity_assumption}, and \ref{Lipschitz_assumption} hold true, for all $t\geq t_0$ when the sequence $\|\bbD_{t-1}^{-1/2}\bbg_{t}\|$ satisfies 
\begin{equation}\label{interval_condition}
	\frac{\sqrt{\eta_{t}}(1- \sqrt{\eta_t}) }{ \eps^2\Gamma_2} 
	\leq\left\|\bbD_{t-1}^{-1/2}\bbg_{t}\right\|
	< \frac{1- \sqrt{\eta_t}}{ \eps^2\Gamma_2}\ ,
\end{equation}
the rate of convergence is quadratic in the sense that
\begin{equation}\label{super_linear_claim_theorem}
{\left\|\bbD_{t}^{-1/2}\bbg_{t+1}\right\|} \leq \frac{\eps^2 \Gamma_2}{1-\sqrt{\eta_t}} {\left\|\bbD_{t-1}^{-1/2}\bbg_{t}\right\|}^{2}.
\end{equation} \end{thm}

%
\begin{myproof} Considering the definition of $\eta_{t}$ we can rewrite the result of Lemma \ref{two_phase_convergence_lemma} as
\begin{equation}\label{safa_city}
\left\|\bbD_{t}^{-{1}/{2}}\bbg_{t+1}\right\|   \leq 
	\eta_t \left\|\bbD_{t-1}^{-{1}/{2}}\bbg_{t}\right\| + \eps^2 \Gamma_2 \left\| \bbD_{t-1}^{-{1}/{2}}\bbg_{t}\right\|^2.
\end{equation}
we use this expression to prove that the inequality in \eqref{super_linear_claim_theorem} holds true. To do so rearrange terms in the first inequality in \eqref{interval_condition} and write
\begin{equation}\label{yek_chi_begim_bere}
{\sqrt{\eta_{t}}}\ \leq\ \frac{\eps^2   \Gamma_2}{1- \sqrt{\eta_t}} \left\| \bbD_{t-1}^{-1/2}\bbg_{t}\right\|.
\end{equation}
Multiplying both sides of \eqref{yek_chi_begim_bere} by $\sqrt{\eta_t}\| \bbD_{t-1}^{-1/2}\bbg_{t}\|$ yields
\begin{equation}\label{yek_chi_begim_bere_2}
\eta_t \left\|\bbD_{t-1}^{-{1}/{2}}\bbg_{t}\right\|  
	\leq\frac{\sqrt{\eta_t}\eps^2   \Gamma_2}{1- \sqrt{\eta_t}}\left\| \bbD_{t-1}^{-1/2}\bbg_{t}\right\|^2.
\end{equation}
Substituting $\eta_t\|\bbD_{t-1}^{-\frac{1}{2}}\bbg_{t}\| $ in \eqref{safa_city} for its upper bound in \eqref{yek_chi_begim_bere_2} implies that 
\begin{align}\label{super_linear_phase}
\left\|\bbD_{t}^{-{1}/{2}}\bbg_{t+1}\right\| 
  &\leq \frac{\sqrt{\eta_t}\eps^2   \Gamma_2}{1- \sqrt{\eta_t}} \left\| \bbD_{t-1}^{-1/2}\bbg_{t}\right\|^2
  + \eps^2 \Gamma_2 \left\| \bbD_{t-1}^{-{1}/{2}}\bbg_{t}\right\|^2	\nonumber \\
  &= \frac{\eps^2   \Gamma_2}{1- \sqrt{\eta_t}} \left\| \bbD_{t-1}^{-1/2}\bbg_{t}\right\|^2.
\end{align}
To verify quadratic convergence, it is necessary to prove that the sequence $\| \bbD_{i-1}^{-1/2}\bbg_{i}\|$ of weighted gradient norms is decreasing. For this to be true we must have
\begin{equation}\label{less_than_1}
 \frac{\eps^2   \Gamma_2}{1- \sqrt{\eta_t}} \left\| \bbD_{t-1}^{-1/2}\bbg_{t}\right\|< 1.
\end{equation}
But \eqref{less_than_1} is true because we are looking at a range of gradients that satisfy the second inequality in \eqref{interval_condition}. \end{myproof}

%
As per Theorem \ref{linear_convergence} $\bby_t$ is converging to $\bby^*$ at a rate that is at least linear. Thus, the gradients $\bbg_t$ will be such that at some point in time they satisfy the rightmost inequality in \eqref{interval_condition}. At that point in time, progress towards $\bby^*$ proceeds at a quadratic rate as indicated by \eqref{super_linear_claim_theorem}. This quadratic rate of progress is maintained until the leftmost inequality in \eqref{interval_condition} is satisfied, at which point the linear term in \eqref{convg_rate_lemma_claim} dominates and the convergence rate goes back to linear. We emphasize that the quadratic convergence region is nonempty because  we have $\sqrt{\eta_t}<1$ for all $t\geq t_0$. Furthermore, making $\eps=1$ and $K$ sufficiently large it is possible to reduce $\eta_t$ arbitrarily and make the quadratic convergence region last longer. In practice, this calls for making $K$ large enough so that $\sqrt{\eta_t}$ is close to the desired gradient norm accuracy.

%
\begin{remark}\label{rmk_make_rho_small}\normalfont
Making $\rho^{K+1}$ small reduces the factor in front of the linear term in \eqref{convg_rate_lemma_claim_large_t_eps_1} and makes the quadratic phase longer. This factor, as it follows from the definition in Proposition \ref{symmetric_term_bounds11}, is $\rho^{K+1}= \left[2(1-\delta)/(2(1-\delta)+{\alpha m})\right]^{K+1}$. Thus, other than increasing $K$, we can make $\rho$ small by increasing the product $\alpha m$. That implies making the inverse penalty coefficient $\alpha$ large relative to the smallest Hessian eigenvalue of the local functions $f_i$ [cf. \eqref{local_hessian_eigenvlaue_bounds}]. This is not possible if we want to keep the solution  $\bby^*$ of \eqref{centralized_opt_problem} close to the solution of $\tby^*$ of \eqref{original_optimization_new_notation}. This calls for the use of adaptive rules to decrease the inverse penalty coefficient $\alpha$ as we elaborate in Section \ref{sec:implement}. Further observe that $\rho$ is independent of the condition number $M/m$ of the local objectives. Making $\rho$ small is an algorithmic choice -- controlled by the selection of $\alpha$ and $K$ --, and not a property of the function being minimized.
\end{remark}

%
\begin{remark} \label{rmk_quadratic_obj}\normalfont
For a quadratic objective function $F$, the Lipschitz constant for the Hessian is $L=0$. Then, the optimal choice of stepsize for NN-$K$ is $\eps=1$ as a result of stepsize rule in \eqref{step_size_condition}. Moreover, the constants for the linear and quadratic terms in \eqref{convg_rate_lemma_claim} are $\Gamma_1=\Gamma_2=0$ as it follows from their definitions in \eqref{Gammas_definition}. For quadratic functions we also have that the Hessian of the objective function $\bbH_t=\bbH$ and the block diagonal matrix $\bbD_{t}=\bbD$ are time invariant, which implies that we can rewrite \eqref{convg_rate_lemma_claim} as 
\begin{equation}\label{eqn_qudartic_condition_numner}
   \|\bbD^{-{1}/{2}}\bbg_{t+1}\| \leq \rho^{K+1}\|\bbD^{-{1}/{2}}\bbg_{t}\| .
\end{equation}
We know that when applying Newton's method to quadratic functions we converge in a single step. This property follows from \eqref{eqn_qudartic_condition_numner} because Newton's method is equivalent to NN-$K$ as $K\to\infty$. The expression in \eqref{eqn_qudartic_condition_numner} states that NN-$K$ converges linearly with a constant decrease factor of $\rho^{K+1}$ per iteration. This factor is independent of the condition number of the quadratic function; see Remark \ref{rmk_make_rho_small}. This in contrast with first order methods like DGD that converge with a linear rate that depends on the condition number of the objective. 
\end{remark}


%
\begin{algorithm}[t]{\small
\caption{Network Newton-$K$ method at node $i$}
\label{algo_NNK} 
\begin{algorithmic}[1]
   \STATE \textbf{function}  
           $\bbx_{i}$
          = NN-$K$$\left(\alpha, \bbx_{i},\text{tol} \right)$ 
   \REPEAT
   \STATE $\bbB$ matrix blocks: 
          $\bbB_{ii}=(1-w_{ii})\bbI$ and $\bbB_{ij}=w_{ij}\bbI$
   \STATE $\bbD$ matrix block: 
          $\bbD_{ii}= \alpha \nabla^2 f_{i}(\bbx_{i}) + 2(1-w_{ii})\bbI $
   \STATE Exchange iterates $\bbx_{i}$ with neighbors $j\in \mathcal{N}_i$.
   \STATE Gradient:    
          $\displaystyle{
          \bbg_{i} = (1-w_{ii})\bbx_{i} 
                       - \sum_{j\in \mathcal{N}_i} w_{ij} \bbx_{j}
                       +\alpha \nabla f_{i}(\bbx_{i}).}$  
   \STATE Compute NN-0 descent direction $\bbd_{i}^{(0)}=-\bbD_{ii}^{-1}\bbg_{i}$\\ 
   \FOR  {$k=  0, \ldots, K-1$ } 
      \STATE Exchange elements $\bbd_{i}^{(k)}$ of the NN-$k$ step with neighbors
      \STATE NN-$(k+1)$ step:
             $\displaystyle{  
             \bbd_{i}^{(k+1)} = \bbD_{ii}^{-1}
             			\bigg[ \sum_{j\in \mathcal{N}_i,j=i}\bbB_{ij} \bbd_{j}^{(k)} 
                                    - \bbg_{i}\bigg]}$.        
   \ENDFOR
   \STATE Update local iterate: 
          $\displaystyle{\bbx_{i}=\bbx_{i} +\eps\ \bbd_{i}^{(K)}}$.
   \UNTIL {$\|\bbg_{i}\|<\text{tol}$}
\end{algorithmic}}\end{algorithm}

%
\begin{algorithm}[t]
\caption{Adaptive Network Newton-$K$ method at node $i$}\label{algo_ANN} 
\begin{algorithmic}[1] {\small
\REQUIRE Iterate $\bbx_{i}$. Initial parameter $\alpha$. 
         Flags $s_{ij}=0$. Factor $\eta<1$.
\FOR {$t=0,1,2,\ldots$}
     \STATE Call NN-$K$ function:
          $\bbx_{i}$
          = NN-$K$$\left(\alpha, \bbx_{i}, \text{tol} \right)$ 
   \STATE Set $s_{ii}=1$ and broadcast it to all nodes.
   \STATE Set $s_{ij}=1$ for all nodes $j$ that sent the signal $s_{jj}=1$.
   \IF{$s_{ij}=1$ for all $j=1,\dots,n$}
      \STATE Update penalty parameter $\alpha=\eta\alpha.$
      \STATE Set $s_{ij}=0$ for all $j=1,\dots,n$.
   \ENDIF
\ENDFOR}
\end{algorithmic}\end{algorithm}

%
\section{Implementation issues} \label{sec:implement}

As mentioned in Section \ref{sec:problem}, NN-$K$ does not solve \eqref{original_optimization_problem1} or its equivalent  \eqref{original_optimization_new_notation}, but the penalty version introduced in \eqref{centralized_opt_problem}. 
The optimal solutions of the optimization problems in \eqref{original_optimization_new_notation} and \eqref{centralized_opt_problem} are different and the gap between them is of order $O(\alpha)$, \cite{YuanQing}. This observation implies that by setting a decreasing policy for $\alpha$, or equivalently, an increasing policy for the penalty coefficient $1/\alpha$, the solution of  \eqref{original_optimization_new_notation} approaches the minimizer of  \eqref{centralized_opt_problem}, i.e. $\tby^* \to \bby^*$ for $\alpha\to0$. 

There are various possible alternatives to reduce $\alpha$. Given the penalty method interpretation in Section \ref{sec:problem} it is more natural to consider fixed penalty parameters $\alpha$ that are decreased after detecting convergence to the optimum argument of the function $F(\bby)$ [cf. \eqref{centralized_opt_problem}]. This latter idea is summarized under the name of Adaptive Network Newton-$K$ (ANN-$K$) in Algorithm \ref{algo_ANN} where $\alpha$ is reduced by a given factor $\eta<1$.

Specifically, ANN-K relies on Algorithm \ref{algo_NNK}, which receives an initial iterate $\bbx_i$, a penalty parameter $\alpha$, and a given tolerance $\text{tol}$ (Step 1) and runs the local NN-$K$ iterations in \eqref{update_formula_NN_local} and \eqref{local_descent} for node $i$ until the local gradient norm $\|\bbg_{i}\|$ becomes smaller than $\text{tol}$ (Step 13). The descent iteration in \eqref{update_formula_NN_local} is implemented in Step 12. Implementation of this descent requires access to the NN-$K$ descent direction $\bbd_{i,t}^{(K)}$ which is computed by the loop in steps 7-11. Step 7 initializes the loop by computing the NN-0 step $\bbd_{i,t}^{(0)}=-\bbD_{ii,t}^{-1}\bbg_{i,t}$. The core of the loop is in Step 10 which corresponds to the recursion in \eqref{local_descent}. Step 8 stands for the variable exchange that is necessary to implement Step 7. After $K$ iterations through this loop, the NN-$K$ descent direction $\bbd_{i,t}^{(K)}$ is computed and can be used in Step 12. Both, steps 7 and 10, require access to the local gradient component $\bbg_{i,t}$. This is evaluated in Step 6 after receiving the prerequisite information from neighbors in Step 5. Steps 3 and 4 compute the blocks $\bbB_{ii,t}$, $\bbB_{ij,t}$, and $\bbD_{ii,t}$ that are also necessary in steps 7 and 10. This process is repeated until $\|\bbg_{i}\|<\text{tol}$ (Step 13). Notice however, that if Algorithm \ref{algo_NNK} is called with a variable $\bbx_i$ with $\|\bbg_{i}\|<\text{tol}$ we still run at least one iteration of NN-$K$. 

ANN-K calls Algorithm \ref{algo_NNK} in Step 2 of Algorithm \ref{algo_ANN}. The factor $\alpha$ is subsequently reduced by the factor $\eta<1$ as indicated in Step 6 of Algorithm \ref{algo_ANN} that implements the replacement $\alpha=\eta\alpha$. The rest of Algorithm \ref{algo_ANN} is designed to handle the fact that a small local gradient norm does not necessarily imply a small global gradient norm. To handle this possible mismatch, flag variables $s_{ij}$ are introduced at node $i$ to signal the fact that node $j$ has reached a local gradient $\bbg_j$ with norm $\|\bbg_j\|\leq\text{tol}$. Whenever node $i$ completes a run of Algorithm \ref{algo_NNK} it broadcasts the signal $s_{ii} $ to all other nodes (Step 3) and updates the variables $s_{ij}$ to $s_{ij}=1$ for all the nodes that sent the signals $s_{jj}=1$ while Algorithm \ref{algo_NNK} was executing (Step 4). If all the variables $s_{ij}=1$ (Step 5) it must be that this is true for all nodes and it is thus safe to modify $\alpha$ (Step 6). The flag variables are reset to $s_{ij}=0$ and Algorithm \ref{algo_NNK} is called with the reduced $\alpha$.

As is typical of penalty methods there are tradeoffs on the selection of the initial value of $\alpha$ and the decrease factor $\eta$. Small values of the initial penalty parameter and $\alpha$ and factor $\eta$ results in sequence of approximate problems having solutions $\tby^*$ that are closer to the actual solution $\bby^*$. However, problems with small $\alpha$ may take a large number of iterations to converge if initialized far from the optimum value because the constant $\rho$ approaches 1 when $\alpha$ is small -- as we discussed in Remark \ref{rmk_make_rho_small}. It is therefore better to initialize Algorithm \ref{algo_ANN} with values of $\alpha$ that are not too small and to decrease $\alpha$ by a factor $\eta$ that is not too aggressive. We discuss these tradeoffs in the numerical examples of Section \ref{sec_sims_ann}.


%
\section{Numerical analysis}\label{sec:simulations}
 
We compare the performance of DGD and different versions of network Newton in the minimization of a distributed quadratic objective. The comparison is done in terms of both, number of iterations and number of information exchanges. For each agent $i$ we consider a positive definite diagonal matrix $\bbA_i\in \mbS_{p}^{++}$ and a vector $\bbb_i \in \mbR^{p}$ to define the local objective function $f_i(\bbx)  :=  ({1}/{2}) \bbx^{T} \bbA_i\bbx +\bbb_i^T\bbx$. Therefore, the global cost function $f(\bbx)$ is written as
\begin{equation}\label{example_problem}
  f(\bbx)  := \sum_{i=1}^n  \frac{1}{2} \bbx^{T} \bbA_i \bbx +\bbb_i^{T}\bbx   \ .
\end{equation}
The difficulty of solving \eqref{example_problem} is given by the condition number of the matrices $\bbA_i$. To adjust condition numbers we generate diagonal matrices $\bbA_i$ with random diagonal elements $a_{ii}$. The first $p/2$ diagonal elements $a_{ii}$ are drawn uniformly at random from the discrete set $\{1, 10^{-1},\ldots, 10^{-\xi}\}$ and the next $p/2$ are uniformly and randomly chosen from the set  $\{1,10^1,\ldots, 10^{\xi}\}$. This choice of coefficients yields local matrices $\bbA_i$ with eigenvalues in the interval $[10^{-\xi}, 10^{\xi}]$ and global matrices $\sum_{i=1}^{n}\bbA_i$ with eigenvalues in the interval $[n10^{-\xi}, n10^{\xi}]$. The linear terms $\bbb_i^{T}\bbx$ are added so that the different local functions have different minima. The vectors $\bbb_i$ are chosen uniformly at random from the box $[0,1]^p$. 

For the quadratic objective in \eqref{example_problem} we can compute the optimal argument $\bbx^*$ in closed form. We then evaluate convergence through the relative error that we define as the average normalized squared distance between local vectors $\bbx_i$ and the optimal decision vector $\bbx^*$,
\begin{equation}\label{distance_error}
e_t := \frac{1}{n} \sum_{i=1}^{n}\frac{\| \bbx_{i,t}-\bbx^*\|^2}{\|\bbx^*\|^2}.
\end{equation}
The network connecting the nodes is a $d$-regular cycle where each node is connected to exactly $d$ neighbors and $d$ is assumed even. The graph is generated by creating a cycle and then connecting each node with the $d/2$ nodes that are closest in each direction. The diagonal weights in the matrix $\bbW$ are set to $w_{ii} = 1/2+1/2(d+1)$ and the off diagonal weights to $w_{ij}= 1/2(d+1)$ when $j\in\ccalN_i$. 
%
\begin{figure}[t]
\centering
\includegraphics[width=\linewidth,height=0.63\linewidth]{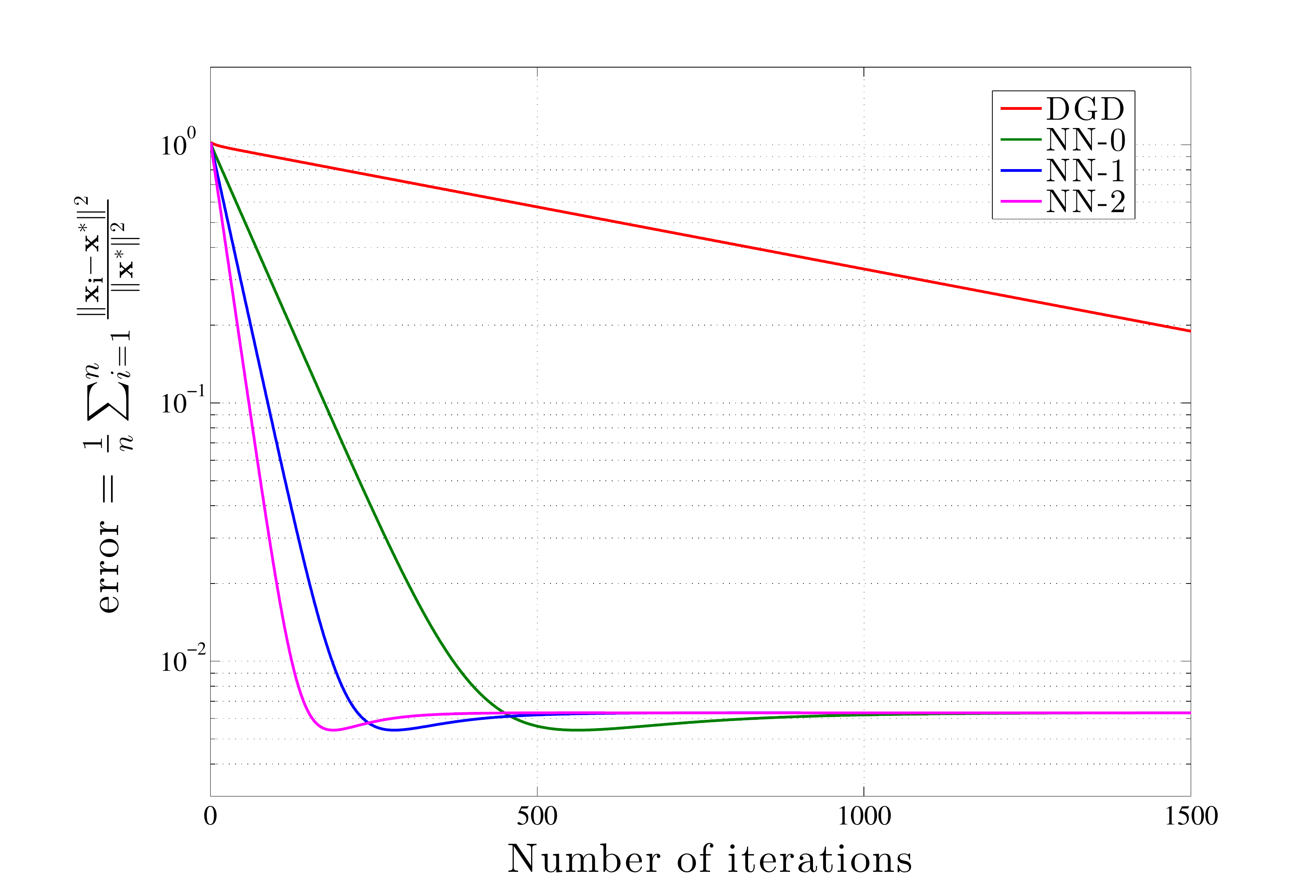}
\caption{Convergence of DGD, NN-0, NN-1, and NN-2 in terms of number of iterations. The network Newton methods converges faster than DGD. Furthermore, the larger $K$ is, the faster NN-$K$ converges.}
\label{fig:iter_illus}
\end{figure}

%
In the subsequent experiments we set the network size to $n=100$, the dimension of the decision vectors to $p=4$, the condition number parameter to $\xi=2$, the penalty coefficient inverse to $\alpha =10^{-2}$, and the network degree to $d=4$. The network Newton step size is set to $\eps=1$, which is always possible when we have quadratic objectives  [cf. Remark \ref{rmk_quadratic_obj}]. Figure \ref{fig:iter_illus} illustrates a sample convergence path for DGD, NN-0, NN-1, and NN-2 by measuring the relative error $e_t$ in \eqref{distance_error} with respect to the number of iterations $t$. As expected for a problem that doesn't have a small condition number -- in this particular instantiation of the function in \eqref{example_problem} the condition number is $95.2$ -- different versions of network Newton are much faster than DGD. E.g., after $t=1.5\times 10^3$ iterations the error associated which DGD iterates is $e_t\approx1.9\times 10^{-1}$. Comparable or better accuracy $e_t<1.9\times 10^{-1}$ is achieved in $t=132$, $t=63$, and $t=43$ iterations for NN-0, NN-1, and NN-2, respectively. 

Further recall that $\alpha$ controls the difference between the actual optimal argument  $\tby^*=[\bbx^*;\ldots;\bbx^*]$ [cf. \eqref{original_optimization_new_notation}] and the argument $\bby^*$ [cf. \eqref{centralized_opt_problem}] to which DGD and network Newton converge. Since we have $\alpha =10^{-2}$ and the difference between these two vectors is of order $O(\alpha)$, we expect the error in \eqref{distance_error} to settle at $e_t\approx10^{-2}$. The error actually settles at $e_t\approx6.3\times 10^{-3}$ and it takes all three versions of network Newton less than $t=400$ iterations to do so. It takes DGD more than $t=10^4$ iterations to reach this value. This relative performance difference decreases if the problem has better conditioning but can be made arbitrarily large by increasing the condition number of the matrix $\sum_{i=1}^{n}\bbA_i$. The number of iterations required for convergence can be further decreased by considering higher order approximations in \eqref{Hessian_approximation_iteration}. The advantages would be misleading because they come at the cost of increasing the number of communications required to approximate the Newton step.

%
\begin{figure}[t]
\centering
\includegraphics[width=\linewidth,height=0.63\linewidth]{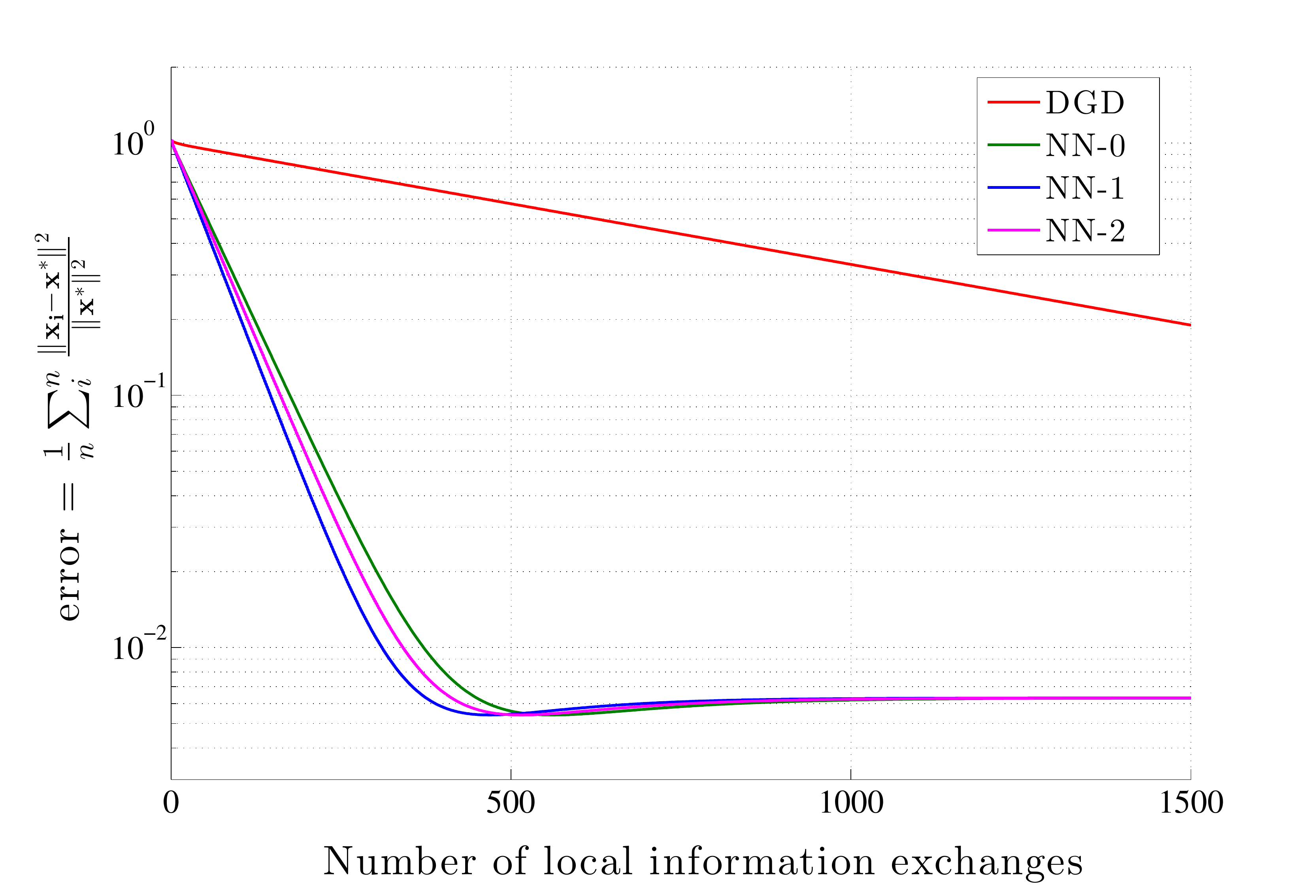}
\caption{Convergence of DGD, NN-0, NN-1, and NN-2 in terms of number of communication exchanges. The NN-$K$ methods retain the advantage over DGD but increasing $K$ may not result in faster convergence. For this particular instance it is actually NN-1 that converges fastest in terms of number of communication exchanges. }
\label{fig:local_exch_illus}
\end{figure}

%
To study this latter effect we consider the relative performance of DGD and different versions of network Newton in terms of the number of local information exchanges. As pointed out in Remark \ref{rmk_comm_cost}, each iteration in NN-$K$ requires a total of $K+1$ information exchanges with each neighbor, as opposed to the single variable exchange required by DGD. After $t$ iterations the number of variable exchanges between each pair of neighbors is $t$ for DGD and $(K+1)t$ for NN-$K$. Thus, we can translate Figure \ref{fig:iter_illus} into a path in terms of number of communications by scaling the time axis by $(K+1)$. The result of this scaling is shown in Figure \ref{fig:local_exch_illus}. The different versions of network Newton retain a significant, albeit smaller, advantage with respect to DGD. Error $e_t<10^{-2}$ is achieved by NN-0, NN-1, and NN-2 after $(K+1)t=3.7\times 10^2$, $(K+1)t=3.1\times 10^2$, and $(K+1)t=3.4\times 10^2$ variable exchanges, respectively. When measured in this metric it is no longer true that increasing $K$ results in faster convergence. For this particular problem instance it is actually NN-1 that converges fastest in terms of number of communication exchanges.

%
\begin{figure}[t]
\centering
\includegraphics[width=\linewidth,height=0.63\linewidth]{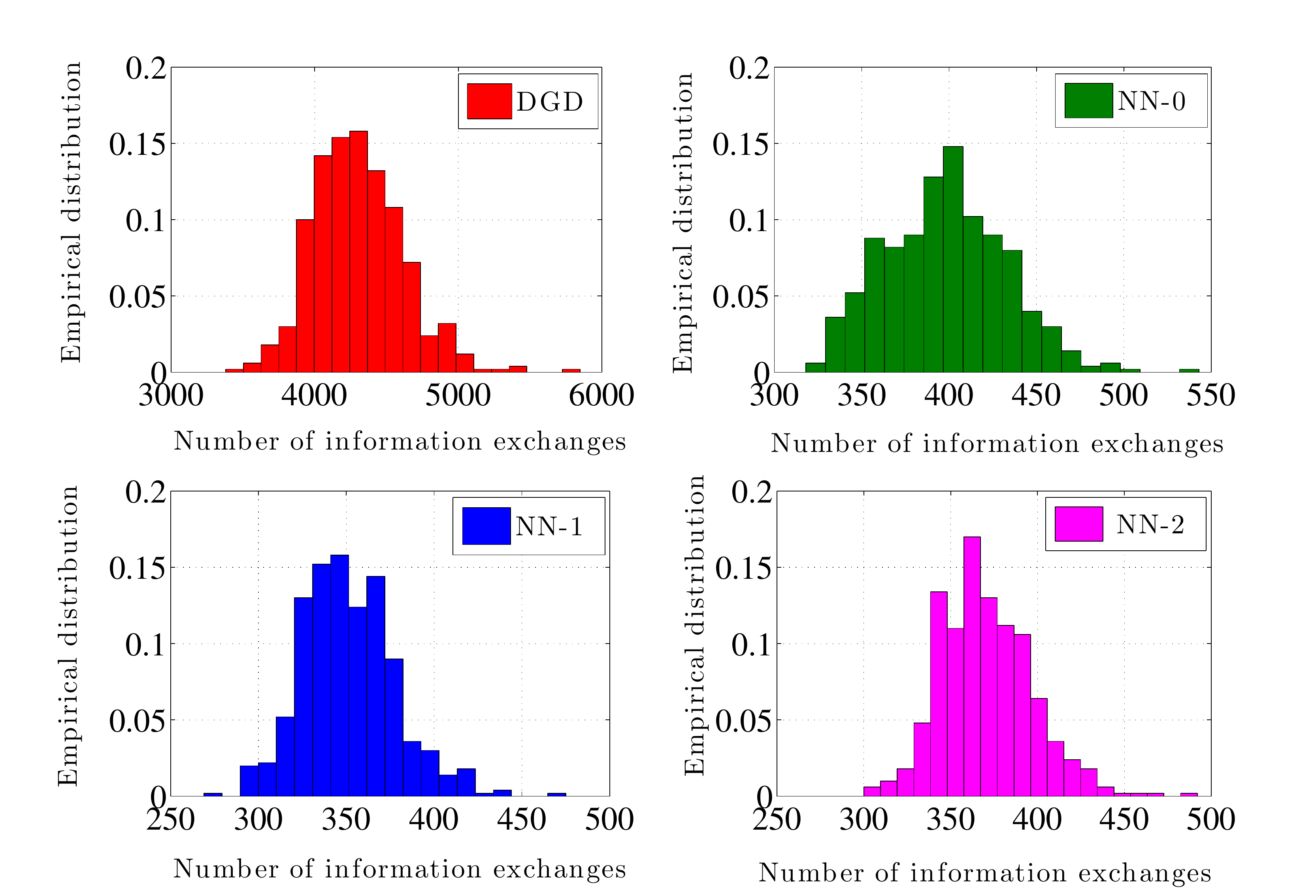}
\caption{Histograms of the number of information exchanges required to achieving accuracy $e_t<10^{-2}$. The qualitative observations made in figures \ref{fig:iter_illus} and \ref{fig:local_exch_illus} hold over a range of random problem realizations.}
\label{fig:local_exch_dist}
\end{figure}

%
For a more comprehensive evaluation we consider $10^3$ different random realizations of \eqref{example_problem} where we also randomize the degree $d$ of the $d$-regular graph that we choose from the even numbers in the set $[2,10]$. The remaining parameters are the same used to generate figures \ref{fig:iter_illus} and \ref{fig:local_exch_illus}. For each joint random realization of network and objective we run DGD, NN-0, NN-1, and NN-2, until achieving error $e_t<10^{-2}$ and record the number of communication exchanges that have elapsed -- which amount to simply $t$ for DGD and $(K+1)t$ for NN. The resulting histograms are shown in Figure \ref{fig:local_exch_dist}. The mean times required to reduce the error to $e_t<10^{-2}$ are $4.3\times 10^3$ for DGD and $4.0\times 10^2$, $3.5\times 10^2$, and $3.7\times 10^2$ for NN-0, NN-1, and NN-2. As in the particular case shown in figures \ref{fig:iter_illus} and \ref{fig:local_exch_illus}, NN-1 performs best in terms of communication exchanges. Observe, however, that the number of communication exchanges required by NN-2 is not much larger and that NN-2 requires less computational effort than NN-1 because the number of iterations $t$ is smaller. 

%
\begin{figure}[t]
\centering
\includegraphics[width=\linewidth,height=0.63\linewidth]{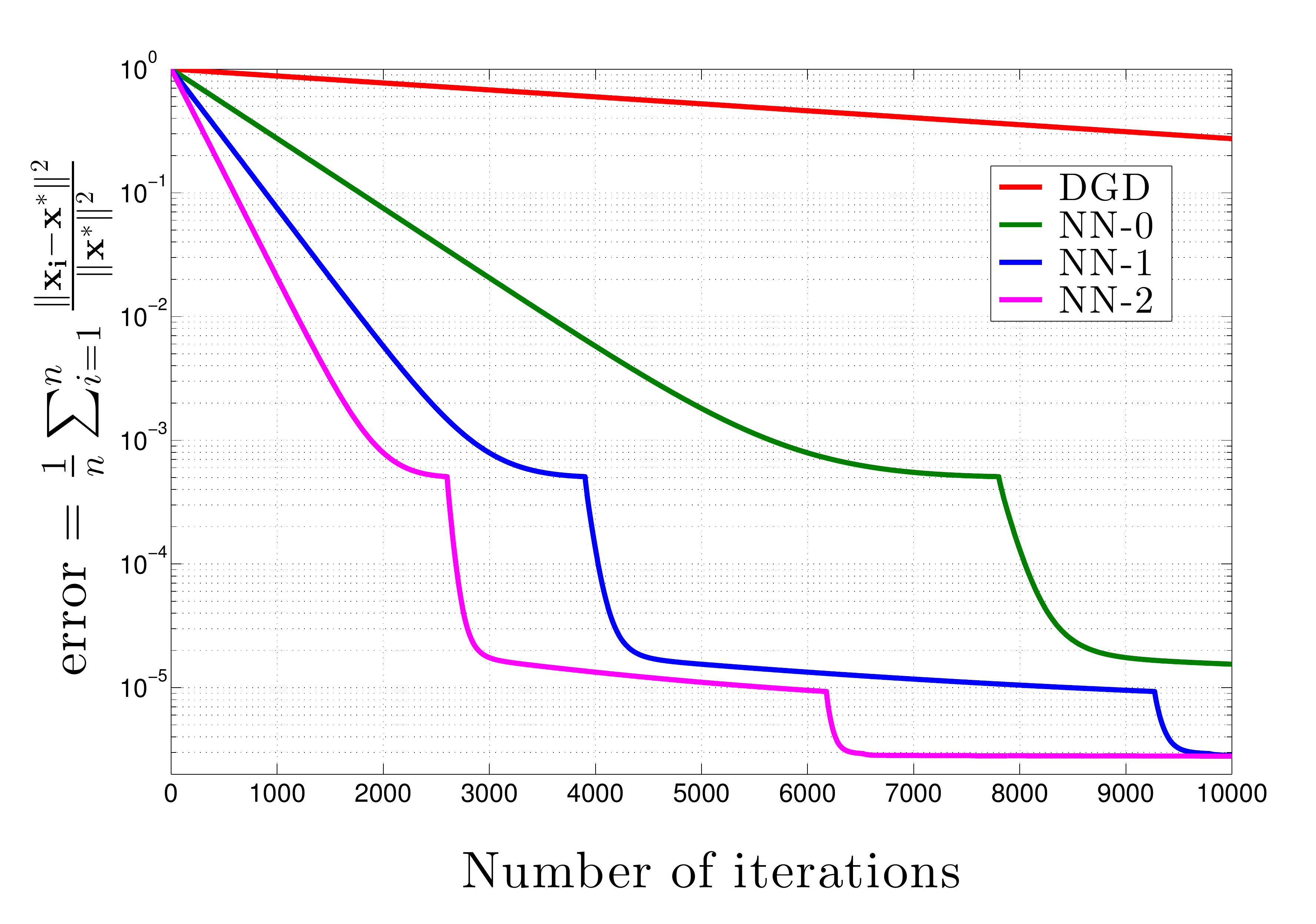}
\caption{{Convergence of adaptive DGD, NN-0, NN-1, and NN-2 for $\alpha_0=10^{-2}$. network Newton methods require less iterations than DGD.}}
\label{fig:4_2}
\end{figure}

%
\subsection{Adaptive Network Newton}\label{sec_sims_ann}

Given that DGD and network Newton are penalty methods it is of interest to consider their behavior when the inverse penalty coefficient $\alpha$ is decreased recursively. The adaptation of $\alpha$ for NN-$K$ is discussed in Section \ref{sec:implement} where it is termed adaptive (A)NN-$K$. The same adaptation strategy is considered here for DGD. The parameter $\alpha$ is kept constant until the local gradient components $\bbg_{i,t}$ become smaller than a given tolerance $\text{tol}$, i.e., until $\|\bbg_{i,t}\|\leq \text{tol}$ for all $i$. When this tolerance is achieved, the parameter $\alpha$ is scaled by a factor $\eta<1$, i.e., $\alpha$ is decreased from its current value to $\eta\alpha$. This requires the use of a signaling method like the one summarized in Algorithm \ref{algo_ANN} for ANN-$K$.

\begin{figure}[t]
\centering
\includegraphics[width=\linewidth,height=0.63\linewidth]{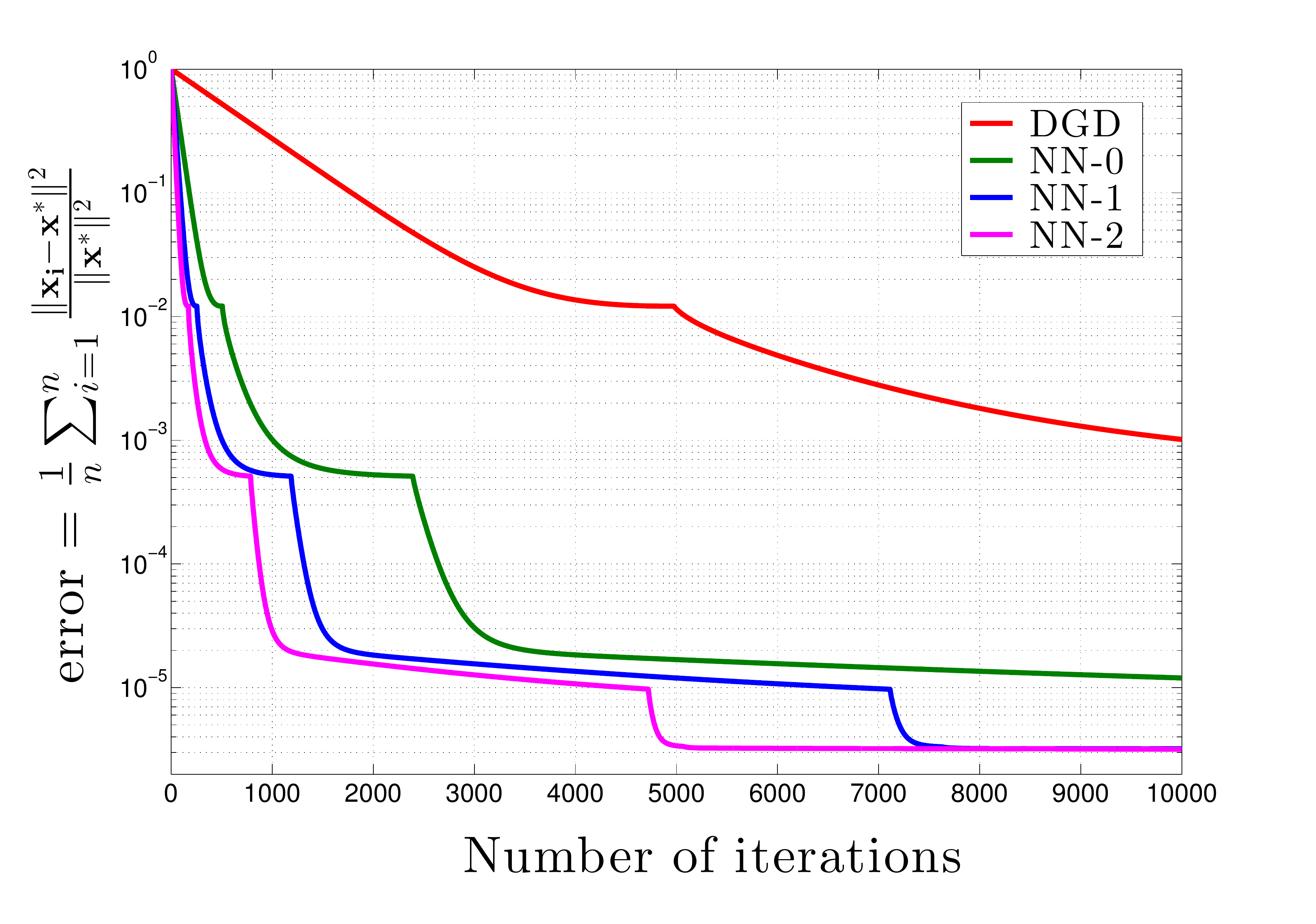}
\caption{{Convergence of Adaptive DGD, NN-0, NN-1, and NN-2 for $\alpha_0=10^{-1}$. ANN methods require less iterations than DGD and convergence of all algorithms are faster relative to the case that $\alpha_0=10^{-2}$.}}
\label{fig:4}
\end{figure}

We consider the objective in \eqref{example_problem} and nodes connected by a $d$-regular cycle. We use the same parameters used to generate figures \ref{fig:iter_illus} and \ref{fig:local_exch_illus}. The adaptive gradient tolerance is set to $\text{tol}=10^{-3}$ and the scaling parameter to $\eta=0.1$. We consider two different scenarios where the initial penalty parameters are $\alpha=\alpha_{0}=10^{-1}$ and $\alpha=\alpha_0=10^{-2}$. The respective error trajectories $e_t$ with respect to the number of iterations are shown in figures \ref{fig:4_2} -- where  $\alpha_0=10^{-2}$ -- and \ref{fig:4} -- where $\alpha_0=10^{-1}$. In each figure we show $e_t$ for adaptive DGD, ANN-0, ANN-1, and ANN-2. Both figures show that the ANN methods outperform adaptive DGD and that larger $K$ reduces the number of iterations that it takes ANN-$K$ to achieve a target error. These results are consistent with the findings summarized in figures \ref{fig:iter_illus}-\ref{fig:local_exch_dist}. 

More interesting conclusions follow from a comparison across figures \ref{fig:4_2} and \ref{fig:4}. We can see that it is better to start with the (larger) value $\alpha=10^{-1}$ even if the method initially converges to a point farther from the actual optimum. This happens because increasing $\alpha$ decreases the constant $\rho=2(1-\delta)/(2(1-\delta)+\alpha m)$. 

%

%
\subsection{Logistic regression}\label{sec_logistic_regression}

For a non-quadratic test we consider the application of network Newton for solving a logistic regression problem. In this problem we are given $q$ training samples that we distribute across $n$ distinct servers. Denote as $q_i$ the number of samples that are assigned to server $i$. Each of the training samples at node $i$ contains a feature vector $\bbu_{il}\in \reals^p$ and a class $v_{il}\in \{-1,1\}$. The goal is to predict the probability  $\Pc{v=1\mid \bbu}$ of having label $v=1$ when given a feature vector $\bbu$ whose class is not known. The logistic regression model assumes that this probability can be computed as $\Pc{v=1\mid \bbu}=1/(1+\exp(-\bbu^T\bbx))$ for a linear classifier $\bbx$ that is computed based on the training samples. It follows from this model that the regularized maximum log likelihood estimate of the classifier $\bbx$ given the training samples $(\bbu_{il},v_{il})$ for $l=1,\ldots,q_i$ and $i=1,\ldots, n$ is given by
\begin{align}\label{eqn_logistic_regrssion_max_likelihood}
   \bbx^*\  =\ &   \argmin_\bbx f(\bbx) \\ \nonumber
         \ :=\ &   \argmin_\bbx   \frac{\lambda}{2} \|\bbx\|^2 
             +               \sum_{i=1}^n \sum_{l=1}^{q_i} 
                                  \log \Big[1+\exp(-v_{il}\bbu_{il}^T\bbx)\Big],
\end{align}
where we defined the function $f(\bbx)$ for future reference. The regularization term $(\lambda/2)\|\bbx\|^2$ is added to reduce overfitting to the training set. 

The optimization problem in \eqref{eqn_logistic_regrssion_max_likelihood} can be written in the form of the optimization problem in \eqref{original_optimization_problem1}. To do so simply define the local objective functions $f_i$ as 
\begin{equation}
   f_i(\bbx) =    \frac{\lambda}{2n} \|\bbx\|^2
                + \sum_{l=1}^{q_i} \log \Big[1+\exp(-v_{il}\bbu_{il}^T\bbx)\Big],
\end{equation}
and observe that given this definition we can write the objective in \eqref{eqn_logistic_regrssion_max_likelihood} as $f(\bbx)=\sum_{i=1}^n f_i(\bbx)$. We can then solve \eqref{eqn_logistic_regrssion_max_likelihood} in a distributed manner using DGD and NN-$K$ methods. 

%
\begin{figure}[t]
\centering
\includegraphics[width=\linewidth,height=0.63\linewidth]{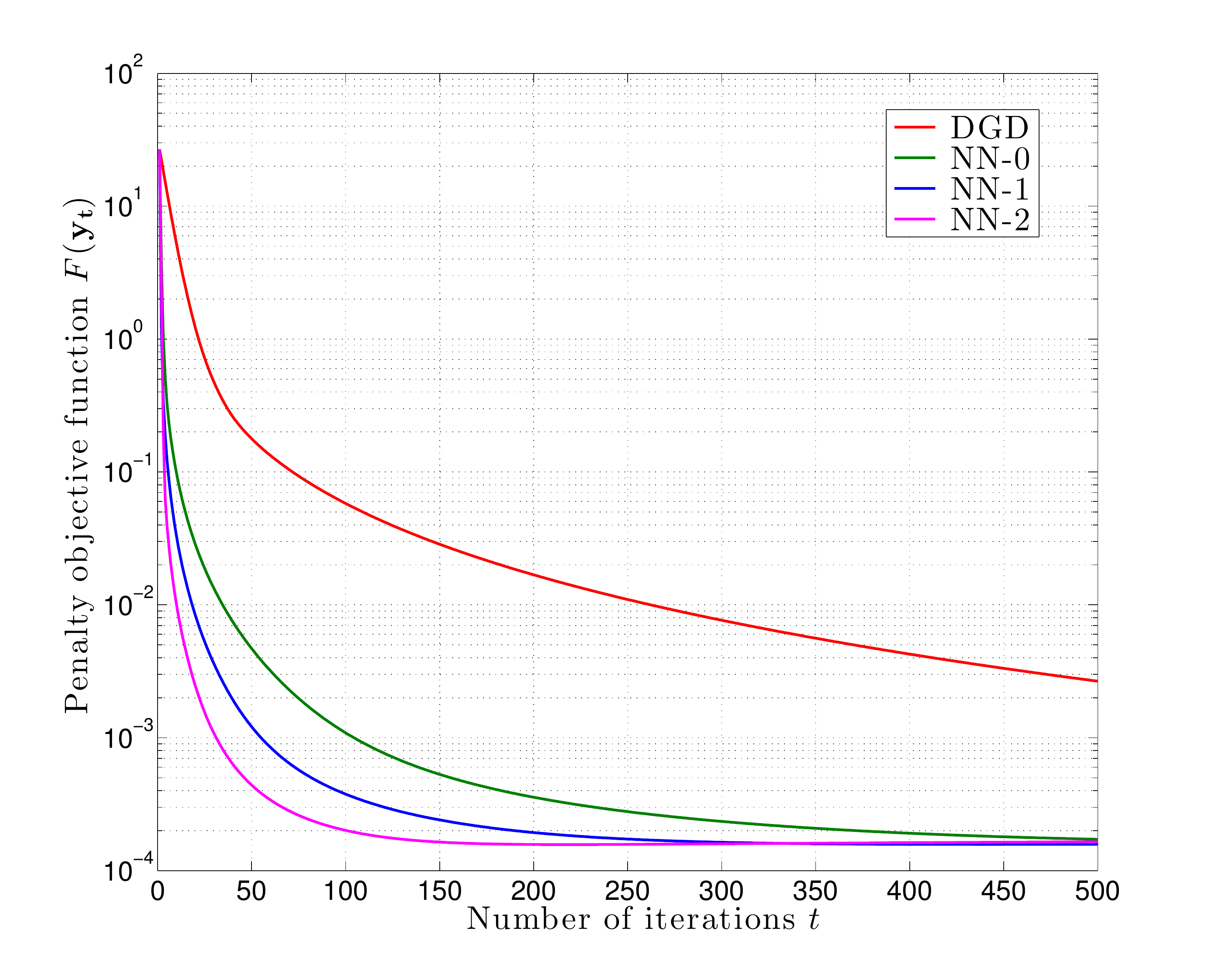}
\caption{{Convergence of DGD, NN-0, NN-1, and NN-2. network Newton methods for a linearly separable logistic regression.}}
\label{fig:5}
\end{figure}

In our experiments we use a synthetic dataset where each component of the feature vector $\bbu_{{il}}$ with label $v_{il}=1$ is generated from a normal distribution with mean $\mu$ and standard deviation $\sigma_+$, while sample points with label $v_{il}=-1$ are generated with mean $-\mu$ and standard deviation $\sigma_-$. The network is a $d$-regular cycle. The diagonal weights in the matrix $\bbW$ are set to $w_{ii} = 1/2+1/2(d+1)$ and the off diagonal weights to $w_{ij}= 1/2(d+1)$ when $j\in\ccalN_i$. We set the feature vector dimension to $p=10$, the number of training samples per node at $q_i=50$, and the regularization parameter to $\lambda=10^{-4}$. The number of nodes is $n=100$ and the degree of the $d$-regular cycle is $d=4$.

We consider first a scenario in which the dataset is linearly separable. To generate a linearly separable dataset the mean is set to $\mu=3$ and the standard deviations to 
$\sigma_+=\sigma_-=1$. Figure \ref{fig:5} illustrates the convergence path of the objective function $F(\bby)$ [cf. \eqref{centralized_opt_problem}] when the penalty parameter is $\alpha=10^{-2}$ and the network Newton step size is $\eps=1$. The reduction in the number of iterations required to achieve convergence is a little more marked than in the quadratic example considered in figures \ref{fig:iter_illus}-\ref{fig:local_exch_dist}. The objective function values $F(\bby_t)$ for NN-0, NN-1 and NN-2 after $t=500$ iterations are below $1.6\times 10^{-4}$, while for DGD the objective function value after the same number of iterations have passed is $F(\bby_t)=2.6\times 10^{-3}$. Conversely, achieving accuracy $F(\bby_t)=2.6\times 10^{-3}$ for NN-0, NN-1, and NN-2 requires $68$, $33$, and $19$ iterations, respectively, while DGD requires $500$ iterations. Observe that for this example NN-2 performs better than NN-1 and NN-0 not only in the number of iterations but also in the number of variable exchanges required to achieve a target accuracy.

\begin{figure}[t]
\centering
\includegraphics[width=\linewidth,height=0.63\linewidth]{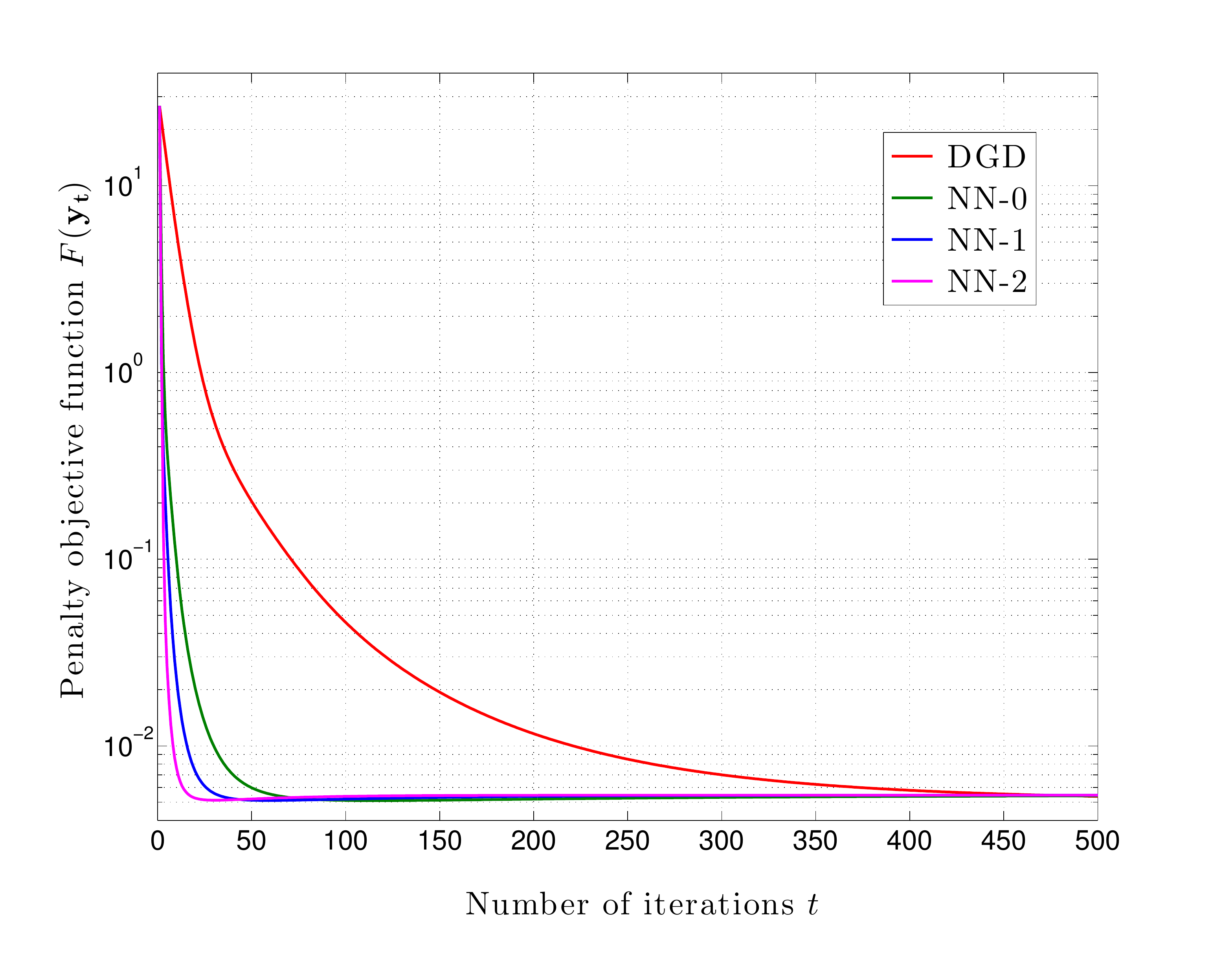}
\caption{{Convergence of DGD, NN-0, NN-1, and NN-2. network Newton methods for a non-linearly separable logistic regression.}}
\label{fig:43}
\end{figure}

We also consider a case in which the dataset is {\it not} linearly separable. To generate this dataset we set the mean to $\mu=2$ and the standard deviations to $\sigma_+=\sigma_-=2$. The penalty parameter is set to $\alpha=10^{-2}$ and the network Newton step size to $\eps=1$. The resulting objective trajectories $F(\bby_t)$ of DGD, NN-0, NN-1, and NN-2 are shown in Figure \ref{fig:43}. The advantages of the network Newton methods relative to DGD are less pronounced but still significant. In this case we also observe that  NN-2 performs best in terms of number of iterations and number of communication exchanges.

%


\section{Conclusions} \label{sec_conclusions}

Network Newton is a decentralized approximation of Newton's method for solving decentralized optimization problems. This paper studied convergence properties and implementation details of this method. Network Newton approximates the Newton direction by truncating a Taylor series expansion of the exact Newton step. This procedure produces a class of algorithms identified by $K$, which is the number of Taylor series terms that network Newton uses for approximating the Newton step. The algorithm is called NN-$K$ when we keep $K$ terms of the Newton step Taylor series. Linear convergence of NN-$K$ is established in a companion paper \cite{NN-part1}. Here, we completed the convergence analysis of NN-$K$ by showing that the sequence of iterates generated by NN-$K$ has a quadratic convergence rate in a specific interval. A quadratic phase exists for all choices of $K$, but this phase can be made arbitrarily large by increasing $K$. The analysis presented here also shows that for the particular case of quadratic objective functions, the convergence rate of NN-$K$ is independent of the condition number of the objective function. This is in contrast to distributed gradient descent methods that require more iterations for problems with larger condition number. Numerical analyses compared the performances of distributed gradient descent and NN-$K$ with different choices of $K$ for minimizing quadratic objectives with large condition numbers as well as the log likelihood objective of a logistic regression problem. In either case we observe that all NN-$K$ methods work faster than distributed gradient descent in terms of number of iterations and number of communications required to achieve convergence. Overall, the theoretical and numerical analyses in this paper prove that NN-$K$ achieves the design goal of accelerating the convergence of distributed gradient descent methods.

We further analyzzed a tradeoff on the selection of a penalty parameter that controls both, the accuracy of the optimal objective computed by network Newton methods and the rate of convergence. We proposed an adaptive version of network Newton (ANN) that achieves exact convergence by executing network Newton with an increasing sequence of penalty coefficients. Numerical analyses of ANN show that it is best to initialize penalty coefficients at moderate values and decrease them through moderate factors.


\begin{appendices}

\section{Proof of Lemma  \ref{two_phase_convergence_lemma}}\label{App_two_phase_convergence_lemma}

{To prove the result in Lemma \ref{two_phase_convergence_lemma}, we first use the Fundamental theorem of Calculus and the Lipschitz continuity of the Hessians $\bbH_t:=\bbH(\bby_t)=\nabla^{2}F(\bby_t)$ to prove the following Lemma.}


\begin{lemma}\label{lemma_help}
Consider the NN-$K$ method as defined in \eqref{diagonal_matrix}-\eqref{update_formula_NN}. If Assumption \ref{Lipschitz_assumption} holds true, then
\begin{equation}\label{lemma_claim_fund_theorem_calculus}
\|\bbg_{t+1}-\bbg_{t}-\bbH_{t}\ (\bby_{t+1}-\bby_{t})\| \leq \frac{\alpha L}{2}\|\bby_{t+1}-\bby_{t}\|^2.
\end{equation}
\end{lemma}

\begin{myproof}
Considering the definitions of the objective function gradient $\bbg_t=\bbg(\bby_t)=\nabla F(\bby_{t})$ and Hessian $\bbH_t=\bbH(\bby_t)=\nabla^2 F(\bby_{t})$, the Fundamental Theorem of calculus implies that 
\begin{equation}\label{eq_1_1}
\bbg_{t+1}=\bbg_{t}+\int_{0}^{1} \bbH(\bby_t + \omega (\bby_{t+1}-\bby_t))\   (\bby_{t+1}-\bby_t) d\omega.
\end{equation}
Adding and subtracting the integral $\int_{0}^{1} \bbH(\bby_t )\   (\bby_{t+1}-\bby_t) d\omega$ to the right hand side of \eqref{eq_1_1} yields 
\begin{align}\label{eq_2_2}
\bbg_{t+1}&=\bbg_{t}+\int_{0}^{1} \bbH(\bby_t )\   (\bby_{t+1}-\bby_t) d\omega \\
&\ + \int_{0}^{1} \left[\bbH(\bby_t + \omega (\bby_{t+1}-\bby_t))-\bbH(\bby_t )\right](\bby_{t+1}-\bby_t) d\omega.\nonumber
\end{align}
Note that $\bbH(\bby_t )\ \!(\bby_{t+1}-\bby_t)$ is not a function of the variable $\omega$ and the first integral in \eqref{eq_2_2} can be simplified as $\int_{0}^{1} \bbH(\bby_t )\   (\bby_{t+1}-\bby_t) d\omega= \bbH(\bby_t )\   (\bby_{t+1}-\bby_t) $. By considering this simplification and regrouping the terms in \eqref{eq_2_2} we obtain
\begin{align}\label{eq_3_3}
&\bbg_{t+1}-\bbg_{t}- \bbH(\bby_t ) (\bby_{t+1}-\bby_t)  =\\ 
&\qquad  \int_{0}^{1} \left[\bbH(\bby_t + \omega (\bby_{t+1}-\bby_t))-\bbH(\bby_t )\right](\bby_{t+1}-\bby_t) d\omega.\nonumber
\end{align}
{We proceed by computing the norm of both sides of \eqref{eq_3_3}. Observing that for any vector $\bba$ the inequality $\int\|\bba\|d\omega\leq \|\int\bba\ \! d\omega\|$ holds and considering that the product of norms is greater than the norm of the respective product, it follows that}
\begin{align}\label{eq_4_4}
&\left\|\bbg_{t+1}-\bbg_{t}- \bbH(\bby_t ) (\bby_{t+1}-\bby_t) \right\| \leq\\ 
&\qquad  \int_{0}^{1} \left\|\bbH(\bby_t + \omega (\bby_{t+1}-\bby_t))-\bbH(\bby_t )\right\|\|\bby_{t+1}-\bby_t\| d\omega.\nonumber
\end{align}
Based on \eqref{H_Lipschitz_claim}, the Hessians $\bbH(\bby_t)$ are Lipschitz continuous with parameter $\alpha L$. Therefore, we can write
\begin{equation}\label{eq_5_5}
\left\|\bbH(\bby_t + \omega (\bby_{t+1}-\bby_t))\!-\!\bbH(\bby_t )\right\| \leq 
\alpha L \omega\|\bby_{t+1}-\bby_t\|.
\end{equation}
Substituting the upper bound in \eqref{eq_5_5} for the term $\left\|\bbH(\bby_t + \omega (\bby_{t+1}-\bby_t))-\bbH(\bby_t )\right\|$ in \eqref{eq_4_4} leads to
\begin{align}\label{eq_6_6}
 \big\|\bbg_{t+1}-\bbg_{t}- \bbH(\bby_t ) &(\bby_{t+1}-\bby_t) \big\| \nonumber\\
    & \leq   \int_{0}^{1} \alpha L \omega\|\bby_{t+1}-\bby_t\|^2 d\omega\nonumber\\
    & =\frac{\alpha L}{2}\|\bby_{t+1}-\bby_t\|^2,
\end{align}
where the equality in \eqref{eq_6_6} is valid since $\int_{0}^{1}\omega\  d\omega=1/2$. The inequality in \eqref{eq_6_6} yields the result in \eqref{lemma_claim_fund_theorem_calculus} considering the notation $\bbH_t=\bbH(\bby_t )$.
\end{myproof}


\textbf{Proof of Lemma \ref{two_phase_convergence_lemma}:} In this proof to simplify the notation we use $\hbH_{t}^{-1}$ to indicate the approximate Hessian inverse $\hbH_{t}^{(K)^{-1}}$. Recall the result in \eqref{lemma_claim_fund_theorem_calculus}. Considering the update formula for NN-$K$ in \eqref{Hessian_approximation_iteration}, the term $\bby_{t+1}-\bby_{t}$ can be substituted by $-\eps \hbH_{t}^{-1} \bbg_{t}$. Making this substitution into \eqref{lemma_claim_fund_theorem_calculus} implies that
\begin{equation}\label{22}
\left\|\bbg_{t+1}-\bbg_{t}+\eps\bbH_{t}\hbH_{t}^{-1} \bbg_{t}\right\| \leq \frac{\eps^2\alpha L}{2}\left\|\hbH_{t}^{-1} \bbg_{t}\right\|^2.
\end{equation}
The definition of matrix norm implies that the norm of product $\bbD_{t}^{-1/2}(\bbg_{t+1}-\bbg_{t}+\eps\bbH_{t}\hbH_{t}^{-1} \bbg_{t})$ is bounded above as
\begin{align}\label{2222}
&\left\|\bbD_{t}^{-{1}/{2}}  \left[\bbg_{t+1}-\bbg_{t}+\eps\bbH_{t}\hbH_{t}^{-1} \bbg_{t}\right]\right\| \leq \nonumber \\
& \qquad \qquad \qquad \left\|\bbD_{t}^{-{1}/{2}} \right\| \left\|\bbg_{t+1}-\bbg_{t}+\eps\bbH_{t}\hbH_{t}^{-1} \bbg_{t}\right\| .
\end{align}
Substituting $\|\bbg_{t+1}-\bbg_{t}+\eps\bbH_{t}\hbH_{t}^{-1} \bbg_{t}\|$ in the right hand side of \eqref{2222} by the upper bound in \eqref{22} leads to 
\begin{align}\label{33}
&\left\|\bbD_{t}^{-{1}/{2}}  \left[\bbg_{t+1}-\bbg_{t}+\eps\bbH_{t}\hbH_{t}^{-1} \bbg_{t}\right]\right\|  \leq  \nonumber \\
&\qquad \qquad \qquad \qquad \quad  \frac{\eps^2 \alpha L}{2}  \left\|\bbD_{t}^{-{1}/{2}} \right\|  \left\|\hbH_{t}^{-1}\bbg_{t}\right\|^2.
\end{align}
Observe that the triangle inequality states that for any vectors $\bba$ and $\bbb $, and a positive constant $C$, if the relation $\|\bba-\bbb\|\leq C$ holds true, then $\|\bba\|\leq \|\bbb\|+C$. By setting $\bba:=\bbD_{t}^{-1/2}\bbg_{t+1}$, $\bbb:=\bbD_{t}^{-1/2} (\bbg_t-\eps\bbH_t\hbH_{t}^{-1})$ and $C:=({\eps^2 \alpha L}/{2})  \|\bbD_{t}^{-{1}/{2}} \|\|\hbH_{t}^{-1}\bbg_{t}\|^2$ and considering the relation in \eqref{33} we obtain that the inequality $\|\bba-\bbb\|\leq C$ is satisfied. Therefore, $\|\bba\|\leq \|\bbb\|+C$ holds true which is equivalent to
\begin{align}\label{444}
\left\|\bbD_{t}^{-{1}/{2}}\bbg_{t+1}\right\| & \leq \left\|\bbD_{t}^{-{1}/{2}}\left[\bbg_{t}-\eps\bbH_{t}\hbH_{t}^{-1} \bbg_{t}\right]\right\| \nonumber \\
&\qquad + \frac{\eps^2 \alpha L}{2}  \left\|\bbD_{t}^{-{1}/{2}}\right\| \left\|\hbH_{t}^{-1}\bbg_{t}\right\|^2.
\end{align}
By rewriting the term $\bbD_{t}^{-1/2} \bbg_t$ as the sum $(1-\eps)(\bbD_{t}^{-1/2} \bbg_t)+\eps(\bbD_{t}^{-1/2} \bbg_t)$ and using the triangle inequality we can update the right hand side of \eqref{444} as
\begin{align}\label{44}
\left\|\bbD_{t}^{-1/2}\bbg_{t+1}\right\| & \leq   (1-\eps)\left\|\bbD_{t}^{-1/2}\bbg_{t}\right\|
	\nonumber \\
	&\qquad
	+ \eps \left\|\bbD_{t}^{-{1}/{2}}\left[\bbI-\bbH_{t}\hbH_{t}^{-1}\right]\bbg_{t}\right\| \nonumber \\
	&\qquad +\frac{\eps^2 \alpha L}{2}  \left\|\bbD_{t}^{-{1}/{2}}\right\| \left\|\hbH_{t}^{-1}\bbg_{t}\right\|^2.
\end{align}
Observe that the Hessian is decomposed as $\bbH_t=\bbD_t-\bbB$ and the approximate Hessian inverse is given by $\hbH_t^{-1} :=     \bbD_t^{-1/2}    \sum_{k=0}^{K} (  \bbD_t^{-1/2}  \bbB   \bbD_t^{-1/2}     )^{k}      \bbD_t^{-1/2}$. Considering these two relations and using the telescopic cancelation we can show that  $\bbI-\bbH_{t}\hbH_{t}^{-1}=(\bbB\bbD_t^{-1})^{K+1}$. This result is studied with more details in Lemma 3 of \cite{NN-part1}. Therefore, we can write
\begin{equation}\label{9909}
\left\|\bbD_{t}^{-1/2}\!\left[\bbI-\bbH_{t}\hbH_{t}^{-1}\right]\!\bbg_{t}\right\|\!
= \!\left\|  \left[\bbD_t^{-1/2}  \bbB\bbD_t^{-1/2} \right]^{K+1}\!\bbD_t^{-1/2}\bbg_{t} \right\|
\end{equation}
{Based on Proposition \ref{symmetric_term_bounds11}, the eigenvalues of the matrix $\bbD_t^{-1/2}  \bbB\bbD_t^{-1/2}$ are bounded by $0$ and $\rho$.} This observation in association with the symmetry of $\bbD_t^{-1/2}  \bbB\bbD_t^{-1/2}$ yields 
\begin{equation}\label{8824}
\left\|\left[\bbD_t^{-1/2}  \bbB\bbD_t^{-1/2}\right]^{K+1}\right\| \leq \rho^{K+1}.
\end{equation}
The simplification in \eqref{9909} and the upper bound in \eqref{8824} guarantee that the norm $\|\bbD_{t}^{-1/2}[\bbI-\bbH_{t}\hbH_{t}^{-1}]\bbg_{t}\|$ is upper bounded by
\begin{equation}\label{7777}
\left\|\bbD_{t}^{-1/2}\left[\bbI-\bbH_{t}\hbH_{t}^{-1}\right]\bbg_{t}\right\|
	 \leq \rho^{K+1} \left\|\bbD_{t}^{-1/2}\bbg_{t}\right\|.
\end{equation}
Substituting the upper bound in \eqref{7777} for the second summand in the right hand side of \eqref{44} implies that
\begin{align}\label{55}
\left\|\bbD_{t}^{-{1}/{2}}\bbg_{t+1}\right\| &\leq   (1-\eps)\left\|\bbD_{t}^{-1/2}\bbg_{t}\right\|
	+ \eps \rho^{K+1} \left\|\bbD_{t}^{-1/2}\bbg_{t}\right\| \nonumber \\
	&\qquad + \frac{\alpha\eps^2  L}{2}  \left\|\bbD_{t}^{-{1}/{2}} \right\| \left\|\hbH_{t}^{-1} \bbg_{t}\right\|^2.
\end{align}
By grouping the first two summands in \eqref{55} and using the inequality $\|\hbH_{t}^{-1} \bbg_{t}\|\leq \|\hbH_{t}^{-1}\|\| \bbg_{t}\|$ we can write 
\begin{align}\label{66}
\left\|\bbD_{t}^{-{1}/{2}}\bbg_{t+1}\right\| 
	&\leq   (1-\eps+ \eps \rho^{K+1}) \left\|\bbD_{t}^{-1/2}\bbg_{t}\right\| \nonumber \\
	&\qquad + \frac{\alpha\eps^2  L}{2}  \left\|\bbD_{t}^{-{1}/{2}} \right\|  \left\|\hbH_{t}^{-1}\right\|^2 \left\| \bbg_{t}\right\|^2.
\end{align}

Now we proceed to find an upper bound for $\|\bbD_{t}^{-1/2}\bbg_t \|$ in terms of $\|\bbD_{t-1}^{-1/2}\bbg_t \|$ to set a recursive inequality for the sequence $\|\bbD_{i-1}^{-1/2}\bbg_i \|$. We first show that the norm of the difference $\|\bbD_{t}^{-1}-\bbD_{t-1}^{-1}\|$ is bounded above as
\begin{align}\label{2200}
\left\|\bbD_{t}^{-1}-\bbD_{t-1}^{-1}\right\| &\leq \left\| \bbD_{t}^{-1}\right\| \left\| \bbD_{t}-\bbD_{t-1}\right\|   \left\| \bbD_{t-1}^{-1}\right\|.
\end{align}
Factoring $\bbD_{t}^{-1}$ and $\bbD_{t-1}^{-1}$ from the left and right sides of $\bbD_{t}^{-1}-\bbD_{t-1}^{-1}$, respectively, follows the relation in \eqref{2200}. Observe that the eigenvalues of the matrices $\bbD_t$ and $\bbD_{t-1}$ are bounded below by $\alpha m+2(1-\Delta)$. Consequently, the eigenvalues of the matrices $\bbD_t^{-1}$ and $\bbD_{t-1}^{-1}$ are bounded above by $1/(\alpha m+2(1-\Delta))$. Therefore, we can update the upper bound in \eqref{2200} as
\begin{equation}\label{2202}
\left\|\bbD_{t}^{-1}-\bbD_{t-1}^{-1}\right\| \leq  \frac{1}{(2(1-\Delta)+\alpha m)^2} \left\| \bbD_{t}-\bbD_{t-1}\right\| .
\end{equation}

The next step is to show that the block diagonal matrices $\bbD_t$ are Lipschitz continuous with parameter $\alpha L$. Notice that the difference $\bbD_t-\bbD_{t-1}$ can be simplified as $\alpha (\bbG_{t}-\bbG_{t-1})$. Moreover, the difference of two consecutive Hessians can be simplified as $\bbH_t-\bbH_{t-1}=\alpha (\bbG_{t}-\bbG_{t-1})$. Therefore, we obtain that $\bbD_t-\bbD_{t-1}=\bbH_t-\bbH_{t-1}$. This observation in association with the Lipschitz continuity of the Hessians with parameter $\alpha L$, i.e., $\|\bbH_t-\bbH_{t-1}\|\leq \alpha L \|\bby_{t}-\bby_{t-1}\|$, implies that 
\begin{equation}\label{D_lip_con}
\left\| \bbD_{t}-\bbD_{t-1}\right\|\leq \alpha L \|\bby_{t}-\bby_{t-1}\|.
\end{equation}
By substituting the upper bound in \eqref{D_lip_con} for the norm $\left\| \bbD_{t}-\bbD_{t-1}\right\|$ in \eqref{2202} we obtain that
\begin{equation}\label{9293}
\left\|\bbD_{t}^{-1}-\bbD_{t-1}^{-1}\right\|
	 \leq \frac{\alpha L}{(2(1-\Delta)+\alpha m)^2}  \left\|\bby_{t}-\bby_{t-1}\right\|.
\end{equation}
Note that the absolute value of the inner product $|\bbg_{t}^{T}(\bbD_{t}^{-1}-\bbD_{t-1}^{-1})\bbg_{t}|$ is bounded above by the product $\|\bbD_{t}^{-1}-\bbD_{t-1}^{-1}\|\|\bbg_t\|^2$.
Considering the upper bound for $\|\bbD_{t}^{-1}-\bbD_{t-1}^{-1}\|$ in \eqref{9293}, the term $|\bbg_{t}^{T}(\bbD_{t}^{-1}-\bbD_{t-1}^{-1})\bbg_{t}|$ is bounded above by
\begin{equation}\label{0912}
\left|\bbg_{t}^{T}(\bbD_{t}^{-1}-\bbD_{t-1}^{-1})\bbg_{t} \right|\leq \frac{\alpha L\  \|\bby_{t}-\bby_{t-1}\| \|\bbg_{t}\|^2}{(2(1-\Delta)+\alpha m)^2} .
\end{equation}
Considering the triangle inequality, and observing the simplifications $|\bbg_t^T \bbD_{t-1}^{-1} \bbg_t|= \|\bbD_{t-1}^{-1/2}\bbg_{t}\|^2$ and $|\bbg_t^T \bbD_{t}^{-1} \bbg_t|= \|\bbD_{t}^{-1/2}\bbg_{t}\|^2$, we can rewrite \eqref{0912} as
\begin{equation}\label{77}
\left\|\bbD_{t}^{-{1}/{2}}\bbg_{t}\right\|^2\leq 
	\left\|\bbD_{t-1}^{-{1}/{2}}\bbg_{t}\right\|^2 
	+\frac{\alpha L\ \|\bby_{t}-\bby_{t-1}\| \|\bbg_{t}\|^2}{(2(1-\Delta)+\alpha m)^2} .
\end{equation}
Observe that for any three constants $a$, $b$ and $c$, if the relation $a^2\leq b^2+c^2$ holds, then the inequality $|a|\leq |b|+|c|$ is valid. By setting $a:=\|\bbD_{t}^{-{1}/{2}}\bbg_{t}\|$, $b:=\|\bbD_{t-1}^{-{1}/{2}}\bbg_{t}\|$, and $c:={\left(\alpha L \|\bby_{t}-\bby_{t-1}\|\right)^{{1}/{2}}\|\bbg_{t}\|}/({2(1-\Delta)+\alpha m})$, we realize that $a^2\leq b^2+c^2$ holds true according to \eqref{77}. Hence, we obtain that the relation $|a|\leq |b|+|c|$ holds and we obtain that 
\begin{equation}\label{88}
\left\|\bbD_{t}^{-{1}/{2}}\bbg_{t}\right\|\leq 
	\left\|\bbD_{t-1}^{-{1}/{2}}\bbg_{t}\right\|
	+\frac{\left(\alpha L \|\bby_{t}-\bby_{t-1}\|\right)^{{1}/{2}}\|\bbg_{t}\|}{2(1-\Delta)+\alpha m}
\end{equation}
Considering the update of NN-$K$ in \eqref{update_formula_NN} we can substitute $\bby_{t}-\bby_{t-1}$ by $-\eps \hbH_{t-1}^{-1} \bbg_{t-1}$. Applying this substitution into \eqref{88} implies that
\begin{equation}\label{99}
\left\|\bbD_{t}^{-{1}/{2}}\bbg_{t}\right\|\leq 
	\left\|\bbD_{t-1}^{-{1}/{2}}\bbg_{t}\right\|
	+\frac{\left[\alpha\eps L \left\|\hbH_{t-1}^{-1} \bbg_{t-1}\right\|\right]^{1/2} \|\bbg_{t}\|}{2(1-\Delta)+\alpha m}.
\end{equation}
If we substitute $\|\bbD_{t}^{-1/2}\bbg_{t}\|$ by the upper bound in \eqref{99} and substitute $\|\hbH_{t-1}^{-1} \bbg_{t-1}\|$ by the upper bound $\|\hbH_{t-1}^{-1}\|\| \bbg_{t-1}\|$, the inequality in \eqref{66} can be written as
\begin{align}\label{100}
&\left\|\bbD_{t}^{-{1}/{2}}\bbg_{t+1}\right\| 
	\leq \left(1-\eps+ \eps \rho^{K+1}\right) \left\|\bbD_{t-1}^{-{1}/{2}}\bbg_{t}\right\| \nonumber \\
	   &\qquad +\frac{ \left(1-\eps+ \eps \rho^{K+1}\right)\left[\alpha\eps L\left\|\hbH_{t-1}^{-1}\right\| \left\|\bbg_{t-1}\right\|\right]^{1/2}}{2(1-\Delta)+\alpha m} \|\bbg_{t}\| \nonumber \\
	   &\qquad + \frac{\alpha\eps^2  L}{2}  \left\|\bbD_{t}^{-{1}/{2}} \right\| \left\|\hbH_{t}^{-1}\right\|^2 \left\| \bbg_{t}\right\|^2.
\end{align}
Due to the fact that for a positive definite matrix the norm of its product by a vector is always larger than its minimum eigenvalue multiplied by the norm of the vector, we can write $\mu_{min}(\bbD_{t-1}^{-1/2})\|\bbg_{t} \| \leq \|\bbD_{t-1}^{-1/2}\bbg_{t}\|$. Rearranging the terms yields
\begin{equation}\label{100202}
\|\bbg_{t} \| \leq \frac{1}{\mu_{min}(\bbD_{t-1}^{-1/2})}\left\|\bbD_{t-1}^{-1/2}\bbg_{t}\right\|. 
\end{equation}
Note that the eigenvalues of the matrix $\bbD_{t-1}$ are upper bounded by $2(1-\delta)+\alpha M$. Hence, $1/\sqrt{(2(1-\delta)+\alpha M)} $ is a lower bound for the eigenvalues of the matrix $\bbD_{t-1}^{-1/2}$. This observation implies that the upper bound in \eqref{100202} can be updated as 
\begin{equation}\label{3388}
\|\bbg_{t}\| \leq (2(1-\delta)+\alpha M)^{1/2} \left\|\bbD_{t-1}^{-{1}/{2}}\bbg_{t}\right\|.
\end{equation}
Substituting $\|\bbg_t\|$ by the upper bound in \eqref{3388} and considering the definition {$\lambda: =1/(2(1-\delta)+\alpha M)$} follows that we can update the right hand side of \eqref{100} as
\begin{align}\label{101}
&\left\|\bbD_{t}^{-{1}/{2}}\bbg_{t+1}\right\| 
         \leq  (1-\eps+ \eps \rho^{K+1}) \left\|\bbD_{t-1}^{-{1}/{2}}\bbg_{t}\right\|  \nonumber \\
         &\quad + \frac{ (1-\eps+ \eps \rho^{K+1})}{(2(1-\Delta)+\alpha m)} \left[\frac{\alpha\eps L\left\|\bbg_{t-1}\right\| \left\|\hbH_{t-1}^{-1}\right\|}{\lambda}\right]^{{1}/{2}}\!\! \left\|\bbD_{t-1}^{-{1}/{2}}\bbg_{t}\right\|
          \nonumber \\
	   & \quad + \frac{\alpha\eps^2  L}{2\lambda}  \left\|\bbD_{t}^{-{1}/{2}} \right\| \left\|\hbH_{t}^{-1}\right\|^2 \left\|\bbD_{t-1}^{-{1}/{2}}\bbg_{t}\right\|^2.
\end{align}
Observe that the norms $\|\hbH_{t}^{-1}\|$ and $\|\hbH_{t-1}^{-1}\|$ are upper bounded by $\Lambda$ which is defined in \eqref{definition_of_lambdas}. Moreover, the norm $\|\bbD_{t}^{-1/2}\|$ is bound above by ${1/(2(1-\Delta)+\alpha m)^{1/2}}$. Substituting these bounds into \eqref{101} results in 
\begin{align}\label{10002}
\left\|\bbD_{t}^{-\frac{1}{2}}\bbg_{t+1}\right\| 
         &\leq  \left(1-\eps+ \eps \rho^{K+1}\right) \left(1+C_1{ \|\bbg_{t-1}\|}^{\frac{1}{2}} \right) \left\|\bbD_{t-1}^{-\frac{1}{2}}\bbg_{t}\right\|   
	\nonumber \\ & \qquad+  \frac{ \alpha\eps^2 L\Lambda^2}{2\lambda {(2(1-\Delta)+\alpha m)}^{\frac{1}{2}}}   \left\| \bbD_{t-1}^{-\frac{1}{2}}\bbg_{t}\right\|^2,
\end{align}
where $C_1$ is defined as
\begin{equation}\label{def_of_C_1}
  C_1 := \left[{\frac{{\alpha \eps L\Lambda}}{\lambda(2(1-\Delta)+\alpha m)^2 }}\right]^{\frac{1}{2}}.
\end{equation}

The next step is to establish an upper bound for ${ \|\bbg_{t-1}\|}^{{1}/{2}} $ in terms of the objective function error $F(\bby_t )-F(\bby^*) $. Observe that the eigenvalues of the Hessian are bounded above by $2(1-\delta)+\alpha M$. This bound in association with the Taylor's expansion of the objective function $F(\bby)$ around $\hby$ leads to 
\begin{equation}\label{taylor_upper_bound}
   F(\bby) \leq\ F(\hby) +\nabla F(\hby)^{T}(\bby-\hby)   
   + {{{2(1-\delta)+\alpha M}\over{2}}\|{\bby - \hby}\|^{2}.}
\end{equation}
According to the definition of $\lambda $ in \eqref{definition_of_lambdas} we can substitute $1/(2(1-\delta)+\alpha M)$ by $\lambda$. Applying this substitution into \eqref{taylor_upper_bound} and minimizing the both sides of \eqref{taylor_upper_bound} with respect to $\bby$ yields
\begin{equation}\label{taylor_upper_bound_2}
   F(\bby^*) \leq\ F(\hby) -\lambda \| \nabla F(\hby) \|^2.
\end{equation}
Since \eqref{taylor_upper_bound_2} holds for any $\hby$, we set $\hby:=\bby_{t-1}$. By rearranging the terms and taking their square roots, we obtain an upper bound for the gradient norm $ \| \nabla F(\bby_{t-1}) \|=\|\bbg_{t-1}\|$ as
\begin{equation}\label{taylor_upper_bound_3}
 \| \bbg_{t-1}\| \leq  \left[\frac{1}{\lambda}  [F(\bby_{t-1}) - F(\bby^*)] \right]^{\frac{1}{2}} .
\end{equation}
The linear convergence of the objective function error implies that $F(\bby_{t-1}) -F(\bby^*) \leq (1-\zeta)^{t-1}(F(\bby_{0}) -F(\bby^*) )$ -- see Theorem \ref{linear_convergence}. Considering this inequality and the relation in \eqref{taylor_upper_bound_3} we can write 
\begin{equation}\label{20202} 
\|\bbg_{t-1}\|^2  \leq \frac{(1-\zeta)^{t-1}}{\lambda}  {\left(F(\bby_{0}) -F(\bby^*) \right)}.
\end{equation}
The upper bound for the squared norm $\|\bbg_{t-1}\|^2$ in \eqref{20202} shows that $\|\bbg_{t-1}\|^{1/2}$ is upper bounded by
\begin{equation}\label{11234}
{\|\bbg_{t-1}\| }^{\frac{1}{2}} \leq   \left[\frac{(1-\zeta)^{t-1}}{\lambda}  {\left(F(\bby_{0}) -F(\bby^*) \right)}\right]^{\frac{1}{4}}.
\end{equation}
By considering the definition of $\Gamma_2$ in \eqref{Gammas_definition} and substituting the upper bound in \eqref{11234} for $\|\bbg_{t-1}\|^{1/2}$, we can update the right hand of \eqref{10002} as
\begin{align}\label{10011111}
\left\|\bbD_{t}^{-\frac{1}{2}}\bbg_{t+1}\right\| 
         &\leq  \left(1-\eps+ \eps \rho^{K+1}\right) \left[1+C_2(1-\zeta)^{\frac{t-1}{4}} \right] \left\|\bbD_{t-1}^{-\frac{1}{2}}\bbg_{t}\right\|   
	\nonumber \\ & \qquad+ \eps^2 \Gamma_2   \| \bbD_{t-1}^{-\frac{1}{2}}\bbg_{t}\|^2,
\end{align}
where $C_2:=C_1[  {(F(\bby_{0}) -F(\bby^*)) /\lambda}]^{1/4}$. Considering the definition of $C_1$ in \eqref{def_of_C_1}, $C_2$ is given by 
\begin{equation}\label{C_2_def}
  C_2 := {\frac{(\alpha \eps L\Lambda)^{\frac{1}{2}} {(F(\bby_{0}) -F(\bby^*)) }^{\frac{1}{4}}}{\lambda^{\frac{3}{4}}(2(1-\Delta)+\alpha m) }}.
\end{equation}
The explicit expression for $C_2$ in \eqref{C_2_def} and the definition of $\Gamma_1$ in \eqref{Gammas_definition} show that $C_2=\Gamma_1$. This observation in association with \eqref{10011111} leads to the claim in \eqref{convg_rate_lemma_claim}.

\end{appendices}
\bibliographystyle{IEEEtran}
  \bibliography{bmc_article}
   \end{document}